\documentclass[aos,MSNbibl,seceqn,citesort,dvips]{arximspdf}
\usepackage{mathrsfs}
\usepackage{graphicx}

\doi{10.1214/11-AOS921}
\volume{40}
\issue{1}
\pubyear{2012}
\firstpage{206}
\lastpage{237}

\makeatletter

\newcommand{\bigmv}{\mid}

\newcommand{\GG}{{\mathbb G}}
\newcommand{\PP}{{\mathbb P}}
\newcommand{\RR}{{\mathbb R}}
\newcommand{\scrB}{{\mathscr B}}
\newcommand{\scrP}{{\mathscr P}}
\newcommand{\score}{\dot{\ell}}
\newcommand{\effscore}{\tilde{\ell}}
\newcommand{\effFI}{{\tilde{I}}}
\newcommand{\effDelta}{{\tilde{\Delta}}}

\newcommand{\ep}{\varepsilon}
\newcommand{\htht}{{\hat{\theta}}}
\newcommand{\Tht}{{\Theta}}

\newcommand{\ft}[2]{\frac{#1}{#2}}
\newcommand{\convprob}[1]{\stackrel{#1}{\longrightarrow}}

\newtheorem{Thm}{Theorem}[section]
\newtheorem{lemma}[Thm]{Lemma}
\newtheorem{corollary}[Thm]{Corollary}

\makeatother

\begin{document}
\begin{frontmatter}

\title{The semiparametric Bernstein--von Mises theorem}
\runtitle{The semiparametric Bernstein--von Mises theorem}

\begin{aug}
\author[A]{\fnms{P. J.} \snm{Bickel}\ead[label=e1]{bickel@stat.berkeley.edu}}
\and
\author[B]{\fnms{B. J. K.} \snm{Kleijn}\corref{}\thanksref{t1}\ead[label=e2]{b.kleijn@uva.nl}}
\runauthor{P. J. Bickel and B. J. K. Kleijn}
\affiliation{University of California, Berkeley and
University of Amsterdam}
\dedicated{Dedicated to the memory of David A. Freedman}
\address[A]{Department of Statistics\\
University of California, Berkeley\\
367 Evans Hall\\
Berkeley, California 94710-3860\\
USA\\
\printead{e1}}
\address[B]{Korteweg-de Vries Institute\\
University of Amsterdam\\
P.O. Box 94248\\
1090 GE, Amsterdam\\
The Netherlands\\
\printead{e2}} %adresu isvedimo komanda gale!
\end{aug}

\thankstext{t1}{Supported by a VENI-grant,
Netherlands Organisation for Scientific Research (NWO).}

% HISTORY:
\received{\smonth{10} \syear{2010}}
\revised{\smonth{8} \syear{2011}}

% ABSTRACT
%
\begin{abstract}
In a smooth semiparametric estimation problem, the marginal posterior
for the parameter of interest is expected to be asymptotically normal
and satisfy frequentist criteria of optimality if the model is endowed
with a suitable prior. It is shown that, under certain straightforward
and interpretable conditions, the assertion of Le Cam's acclaimed, but
strictly parametric, Bernstein--von Mises theorem [\textit{Univ.
California Publ. Statist.} \textbf{1} (1953) 277--329] holds in the
semiparametric situation as well. As a consequence, Bayesian
point-estimators achieve efficiency, for example, in the sense of
H\'ajek's convolution theorem [\textit{Z. Wahrsch. Verw. Gebiete}
\textbf{14} (1970) 323--330]. The model is required to satisfy
differentiability and metric entropy conditions, while the nuisance
prior must assign nonzero mass to certain Kullback--Leibler
neighborhoods [Ghosal, Ghosh and van der Vaart
\textit{Ann. Statist.} \textbf{28}
(2000) 500--531]. In addition, the marginal posterior is required to
converge at parametric rate, which appears to be the most stringent
condition in examples. The results are applied to estimation of the
linear coefficient in partial linear regression, with a~Gaussian prior
on a smoothness class for the nuisance.
\end{abstract}

% KEYWORDS
%

\end{frontmatter}

%%%%%%%%%%%%%%%%%%%%%%%%%%%%%%%%%%%%%%%%%%%%%%%%%%%%%%%%%%%%%%%%%%%%%%%%%%%%%%%%%%
%s1 #&#
\section{Introduction}
\label{secintro}

The concept of efficiency has its origin in Fisher's 1920s
claim of asymptotic optimality of the maximum-likelihood estimator in
differentiable parametric models (Fisher~\cite{Fisher59}). In
1930s and 1940s, Fisher's ideas on optimality in
differentiable models were sharpened and elaborated upon (see, e.g.,
Cram\'er~\cite{Cramer46}), until Hodges's 1951 discovery of a
superefficient estimator indicated that a comprehensive understanding
of optimality in differentiable estimation problems remained elusive.
Further consideration directed attention to the
property of \textit{regularity} to delimit the class of estimators over
which optimality is achieved. H\'ajek's convolution theorem (H\'ajek
\cite{Hajek70}) implies that within the class of regular
estimates, asymptotic variance is lower-bounded by the Cram\'er--Rao
bound in the limit experiment~\cite{LeCam72}. The asymptotic
minimax theorem (H\'ajek~\cite{Hajek72}) underlines the
central role of the concept of regularity. An estimator that
is optimal among regular estimates is called \textit{best-regular};
in a Hellinger differentiable model, an estimator $(\hat{{\theta}}_n)$
for ${\theta}$ is best-regular \textit{if and only if} it is asymptotically
linear, that is, for all ${\theta}$ in the model,
%
%e1.1 #&#
%
\begin{equation}
\label{eqaslin}
\sqrt{n}(\hat{{\theta}}_n-{\theta}) = \frac{1}{\sqrt{n}}\sum_{i=1}^n
I_{\theta}^{-1}\score_{\theta}(X_i) + o_{P_{\theta}}(1),
\end{equation}
where $\score_{\theta}$ is the score for ${\theta}$ and $I_{\theta
}$ the
corresponding Fisher information.
To address the question of efficiency in smooth parametric models
from a Bayesian perspective, we turn to the Bernstein--von Mises
theorem. In the literature many different versions of the theorem
exist, varying both in (stringency of) conditions and (strength or)
form of the assertion. Following Le Cam and Yang~\cite{LeCam90}
(see also van der Vaart~\cite{vdVaart98}), we state the theorem
as follows. (For later reference, define a prior to be \textit{thick}
at ${\theta}_0$, if it has a~Lebesgue density that is continuous and
strictly positive at ${\theta}_0$.)
%
%th1.1 #&#
%
\begin{Thm}[(Bernstein--von Mises, parametric)]
\label{thmparaBvM}
Assume that $\Tht\subset\RR^k$ is open and that the model
$\scrP=\{P_{\theta}\dvtx{\theta}\in\Tht\}$ is identifiable and dominated.
Suppose $X_1, X_2, \ldots$ forms an i.i.d. sample from
$P_{{\theta}_0}$
for some ${\theta}_0\in\Tht$. Assume that the model is locally asymptotically
normal at ${\theta}_0$ with nonsingular Fisher information
$I_{{\theta}_0}$.
Furthermore, suppose that:
\begin{longlist}[(ii)]
\item[(i)] the prior $\Pi_\Tht$ is thick at ${\theta}_0$;
\item[(ii)] for every $\ep>0$, there exists a
test sequence $(\phi_n)$ such that
\[
P_{{\theta}_0}^n\phi_n\rightarrow0,\qquad \sup_{\|{\theta}-{\theta
}_0\|>\ep}
P_{{\theta}}^n(1-\phi_n)\rightarrow0.
\]
\end{longlist}
Then the posterior distributions converge in total variation,
\[
\sup_{B}\bigl| \Pi( {\theta}\in B\bigmv X_1,\ldots
,X_n)
- N_{{\htht_n},(nI_{{\theta}_0})^{-1}}(B) \bigr| \rightarrow0
\]
in $P_{{\theta}_0}$-probability, where $(\htht_n)$ denotes any best-regular
estimator sequence.
\end{Thm}

For a proof, the reader is referred to~\cite{LeCam90,vdVaart98} (or to
Kleijn and van der Vaart~\cite{Kleijn07}, for a proof under model
misspecification that has a lot in common with the proof of
Theorem~\ref{thmpan} below).

Neither the frequentist theory on asymptotic optimality nor
Theorem~\ref{thmparaBvM} generalize fully to nonparametric estimation
problems. Examples of the failure of the Bernstein--von Mises\vadjust{\goodbreak}
limit in infinite-dimensional problems (with regard to the
\textit{full} parameter) can be found
in Freedman~\cite{Freedman99}. Freedman initiated a discussion
concerning the merits of Bayesian methods in nonparametric problems
as early as 1963, showing that even with a natural and seemingly
innocuous choice of the nonparametric prior, posterior inconsistency
may result~\cite{Freedman63}. This warning against instances of
inconsistency due to ill-advised nonparametric priors was reiterated
in the literature many times over, for example, in Cox
\cite{Cox93} and in Diaconis and Freedman
\cite{Diaconis86,Diaconis98}.
However, general conditions for Bayesian consistency were formulated
by Schwartz as early as 1965~\cite{Schwartz65}; positive results on
posterior rates of convergence in the same spirit were obtained in
Ghosal, Ghosh and van der Vaart~\cite{Ghosal00} (see also,
Shen and Wasserman~\cite{Shen01}). The combined
message of negative and positive results appears to be that
the choice of a nonparametric prior is a sensitive one that leaves
room for unintended consequences unless due care is taken.

This lesson must also be taken seriously when one asks the question
whether the posterior for the parameter of interest in a semiparametric
estimation problem displays Bernstein--von Mises-type limiting
behavior. Like in the parametric case, we estimate a finite-dimensional
parameter ${\theta}\in\Tht$, but now in a model $\scrP$ that also
leaves room
for an infinite-dimensional nuisance parameter $\eta\in H$. We look
for general sufficient conditions on model and prior such that the
\textit{marginal posterior for the parameter of interest} satisfies
%
%e1.2 #&#
%
\begin{equation}
\label{eqassertBvM}
\sup_{B}\bigl|
\Pi\bigl( \sqrt{n}({\theta}-{\theta}_0)\in B\bigmv X_1,\ldots
,X_n\bigr)
- N_{\effDelta_n,\effFI_{{\theta}_0,\eta_0}^{-1}}(B)
\bigr| \rightarrow0
\end{equation}
in $P_{{\theta}_0}$-probability, where
%
%e1.3 #&#
%
\begin{equation}
\label{eqDefDeltaMS}
\effDelta_n = \frac{1}{\sqrt{n}}\sum_{i=1}^n
\effFI_{{\theta}_0,\eta_0}^{-1}\effscore_{{\theta}_0,\eta_0}(X_i).
\end{equation}
Here $\effscore_{{\theta},\eta}$ denotes the efficient score
function and
$\effFI_{{\theta},\eta}$ the efficient Fisher information [assumed
to be nonsingular at $({\theta}_0,\eta_0)$]. The sequence $\effDelta_n$
also features on the r.h.s. of the semiparametric version of
(\ref{eqaslin}) (see Lemma 25.23 in~\cite{vdVaart98}). Assertion
(\ref{eqassertBvM}) often implies efficiency of point-estimators
like the posterior median, mode or mean (a first condition being that
the estimate is a functional on $\RR$, continuous in total-variation
\cite{vdVaart98,Kleijn03}) and
always leads to asymptotic identification of credible regions with
efficient confidence regions. To illustrate, if $C$ is a credible set in
$\Tht$,~(\ref{eqassertBvM}) guarantees that posterior coverage
and coverage under the limiting normal for $C$ are (close to) equal.
Because the limiting normals are \textit{also} the asymptotic sampling
distributions for efficient point-estimators,~(\ref{eqassertBvM})
enables interpretation of credible sets as asymptotically efficient
confidence regions. From a practical point of view, the latter
conclusion has an important implication: whereas it can be
hard to compute optimal semiparametric confidence regions
directly, simulation of a large sample from the marginal posterior
(e.g., by MCMC techniques; see Robert~\cite{Robert01})
is sometimes comparatively straightforward.

Instances of the Bernstein--von Mises limit have been studied in
various semiparametric models: several papers have provided studies of
asymptotic normality of posterior distributions for models from
survival analysis. Particularly, Kim and Lee~\cite{Kim04} show that the
\textit{infinite-dimensional} posterior for the cumulative hazard
function under right-censoring converges at rate~$n^{-1/2}$ to a
Gaussian centered at the Aalen--Nelson estimator for a class of
neutral-to-the-right process priors. In Kim~\cite{Kim06}, the posterior
for the baseline cumulative hazard function and regression coefficients
in Cox's proportional hazard model are considered with similar priors.
Castillo~\cite{Castillo08} considers marginal posteriors in Cox's
proportional hazards model and Stein's symmetric location problem from
a unified point of view. A general approach has been given in Shen
\cite{Shen02}, but his conditions may prove somewhat hard to verify in
examples. Cheng and Kosorok~\cite{Cheng08} give a general perspective
too, proving weak convergence of the posterior under sufficient
conditions. Rivoirard and Rousseau~\cite{Rivoirard09} prove a version
for linear functionals over the model, using a~class of nonparametric
priors based on infinite-dimensional exponential families. Boucheron
and Gassiat~\cite{Boucheron09} consider the Bernstein--von Mises
theorem for families of discrete distributions. Johnstone
\cite{Johnstone10} studies various marginal posteriors in the Gaussian
sequence model.
%
%This paper is organised as follows: in sections~\ref{secpert}--
%we discuss the proof of our main result (theorem~\ref{thmsbvmone})
%in three stages and combine them. Section~\ref{secpert} details
%convergence of the nuisance posterior when the parameter of interest
%lies in a $n^{-1/2}$-neighbourhood around its true value. In
%section~\ref{secilan}, we consider a LAN-expansion of the integral
%of the likelihood, used in section~\ref{secpan} to prove asymptotic
%normality of the marginal posterior for the parameter of interest. In
%section~\ref{secmarg} we discuss the asymptotic tail-condition
%for the marginal posterior. In section~\ref{secmainres}, we
%give an overview of the proof and state the main result. We apply
%theorem~\ref{thmsbvmone} in section~\ref{secplr} to the estimation
%of the linear coefficient in the partial linear regression model.

\subsection*{Notation and conventions}

The (frequentist) true distribution of the data is denoted $P_0$ and
assumed to lie in $\scrP$, so that there exist ${\theta}_0\in\Tht$,
$\eta_0\in H$ such that $P_0=P_{{\theta}_0,\eta_0}$. We
localize ${\theta}$ by introducing $h =
\sqrt{n}({\theta}-{\theta}_0)$ with inverse
${\theta}_n(h)={\theta}_0+n^{-1/2}h$. The expectation of a random
variable $f$ with respect to a probability measure $P$ is denoted $Pf$;
the sample average\break of~$g(X)$ is~denoted $\PP_ng(X) =
(1/n)\sum_{i=1}^ng(X_i)$ and $\GG_ng(X) = n^{1/2}(\PP_ng(X)-Pg(X))$ (for
other conventions and nomenclature customary in empirical process
theory, see~\cite{vdVaart96}). If $h_n$ is stochastic,
$P_{{\theta}_n(h_n),\eta}^nf$ denotes the integral $\int f(\omega)
(dP_{{\theta}_n(h_n(\omega )),\eta}^n/ dP_0^n) (\omega) \,dP_0^n(\omega)$.
The Hellinger distance between~$P$\break and~$P'$ is denoted $H(P,P')$ and
induces a metric $d_H$ on the space of nuisance parameters $H$ by
$d_H(\eta,\eta')= H(P_{{\theta}_0,\eta},P_{{\theta}_0,\eta'})$, for all
$\eta,\eta '\in H$. We endow the model with the Borel $\sigma$-algebra
generated by the Hellinger topology and refer to~\cite{Ghosal00}
regarding issues of measurability.

%%%%%%%%%%%%%%%%%%%%%%%%%%%%%%%%%%%%%%%%%%%%%%%%%%%%%%%%%%%%%%%%%%%%%%%%%%%%%%%%%%

%s2 #&#
\section{Main results}
\label{secmainres}

Consider estimation of a functional $\theta\dvtx\scrP
\rightarrow\RR^k$ on a~dominated nonparametric model $\scrP$ with
metric $g$, based on a sample $X_1,X_2,\ldots,$ i.i.d.
according to $P_0\in\scrP$. We introduce
a prior $\Pi$ on $\scrP$ and consider the subsequent sequence of
posteriors,
%
%e2.1 #&#
%
\begin{equation}
\label{eqposterior}
\Pi( A\bigmv X_1,\ldots,X_n)=
{ {\int_A\prod_{i=1}^n p(X_i) \,d\Pi(P)}} \Big/
{ {\int_\scrP\prod_{i=1}^n p(X_i) \,d\Pi(P)}},
\end{equation}
where $A$ is any measurable model subset. Typically, optimal
(e.g., minimax) nonparametric posterior rates of convergence
\cite{Ghosal00} are powers of $n$ (possibly modified by a
slowly varying\vspace*{1pt} function) that converge to zero more
slowly than the parametric $n^{-1/2}$-rate. Estimators for
${\theta}$ may be derived by ``plugging in'' a nonparametric
estimate [cf. $\hat{{\theta}}={\theta}(\hat{P})$], but optimality
in rate
or asymptotic variance cannot be expected to obtain generically
in this way. This does not preclude efficient estimation of
real-valued aspects of $P_0$: parametrize the model in
terms of a finite-dimensional \textit{parameter of interest} ${\theta}\in
\Tht$
and a \textit{nuisance parameter} $\eta\in H$ where $\Tht$ is open
in $\RR^k$ and $(H,d_H)$ an infinite-dimensional metric space:
$\scrP=\{ P_{{\theta},\eta} \dvtx {\theta}\in\Tht,\eta\in H \}$. Assuming
identifiability, there exist unique ${\theta}_0\in\Tht$, $\eta_0\in
H$ such
that $P_0=P_{{\theta}_0,\eta_0}$. Assuming measurability of the map
$({\theta},\eta)\mapsto P_{{\theta},\eta}$, we place a product prior
$\Pi_\Tht\times\Pi_H$ on $\Tht\times H$ to define a prior on
$\scrP$. Parametric rates for the marginal posterior of ${\theta}$
are achievable because it is possible for contraction of the
full posterior to occur anisotropically, that is, at rate $n^{-1/2}$
along the ${\theta}$-direction, but at a slower, nonparametric rate
$(\rho_n)$ along the $\eta$-directions.

%s2.1 #&#
\subsection{Method of proof}
\label{submethod}

The proof of~(\ref{eqassertBvM}) will consist of three steps:
in Section~\ref{secpert}, we show that the posterior concentrates
its mass around so-called \textit{least-favorable submodels} (see
Stein~\cite{Stein56} and~\cite{Bickel98,vdVaart98}). In
the second step (see Section~\ref{secilan}), we show that this
implies local asymptotic normality (LAN) for integrals of the likelihood
over $H$, with the efficient score determining the expansion. In
Section~\ref{secpan}, it is shown that these LAN integrals induce
asymptotic normality of the marginal posterior, analogous to the
way local asymptotic normality of parametric likelihoods
induces the parametric Bernstein--von Mises theorem.

To see why asymptotic accumulation of posterior mass occurs
around so-called least-favorable submodels, a crude argument departs
from the observation that, according to~(\ref{eqposterior}), posterior
concentration occurs in regions of the model with relatively high
(log-)likelihood (barring inhomogeneities of the prior). Asymptotically,
such regions are characterized by close-to-minimal Kullback--Leibler
divergence with respect to $P_0$. To exploit this, let us assume that
for each ${\theta}$ in a neighborhood $U_0$ of ${\theta}_0$, there exists
a unique minimizer $\eta^*({\theta})$ of the Kullback--Leibler divergence,
%
%e2.2 #&#
%
\begin{equation}
\label{eqminKL}
-P_0\log\frac{p_{{\theta},\eta^*({\theta})}}{p_{{\theta}_0,\eta_0}}
= \inf_{\eta\in H} \biggl(-P_0\log\frac{p_{{\theta},\eta
}}{p_{{\theta}_0,\eta
_0}}\biggr)
\end{equation}
giving rise to a submodel $\scrP^*=\{P^*_{\theta}=P_{{\theta},\eta
^*({\theta})}\dvtx
{\theta}\in U_0\}$. As is well known~\cite{Severini92}, if $\scrP^*$ is
smooth it constitutes a least-favorable submodel and scores along~$\scrP^*$
are efficient. [In subsequent sections it is not required
that $\scrP^*$ is defined by~(\ref{eqminKL}), only that $\scrP^*$
is least-favorable.] Neighborhoods of $\scrP^*$ are described
with Hellinger balls in $H$ of radius $\rho>0$ around $\eta^*({\theta})$,
for all ${\theta}\in U_0$,
%
%e2.3 #&#
%
\begin{equation}
\label{eqdefD}
D({\theta},\rho)=\{ \eta\in H \dvtx d_H(\eta,\eta^*({\theta
}))<\rho\}.
\end{equation}
To give a more precise argument for posterior concentration around
$\eta^*({\theta})$, consider the posterior for $\eta$, \textit{given}
${\theta}\in U_0$; unless ${\theta}$ happens to be equal to~${\theta
}_0$, the submodel $\scrP_{\theta}= \{P_{{\theta},\eta}\dvtx\eta\in
H\}$ is misspecified. Kleijn and van der Vaart~\cite{Kleijn06} show
that the misspecified posterior concentrates asymptotically in any
(Hellinger) neighborhood of the point of minimal Kullback--Leibler
divergence with respect to the true distribution of the data. Applied
to $\scrP _{\theta}$, we see that $D({\theta},\rho)$ receives
asymptotic posterior probability one for any $\rho>0$. For posterior
concentration to occur~\cite{Ghosal00,Kleijn06} sufficient prior mass
must be present in certain Kullback--Leibler-type
neighborhoods. In the present context, these neighborhoods can be
defined as
%
%e2.4 #&#
%
\begin{eqnarray}
\label{eqKsets} K_n(\rho,M) &=& \biggl\{ \eta\in H\dvtx P_0\biggl(
\sup_{\|h\|\leq M}-\log\frac{p_{{\theta}_n(h),\eta}}{p_{{\theta
}_0,\eta_0}}
\biggr)\leq\rho^2,\nonumber\\[-8pt]\\[-8pt]
&&\hspace*{37.5pt}P_0\biggl( \sup_{\|h\|\leq
M}-\log\frac{p_{{\theta}_n(h),\eta}}{p_{{\theta }_0,\eta_0}}
\biggr)^2\leq\rho^2 \biggr\}\nonumber
\end{eqnarray}
for $\rho>0$ and $M>0$. If this type of posterior convergence occurs
with an appropriate form of uniformity over the relevant values
of ${\theta}$ (see ``consistency under perturbation,'' Section~\ref{secpert}),
one expects that the nonparametric posterior contracts into Hellinger
neighborhoods of the curve ${\theta}\mapsto({\theta},\eta^*({\theta}))$
(Theorem~\ref{thmpertroc} and Corollary~\ref{corconspert}).

To introduce the second step, consider~(\ref{eqposterior}) with
$A=B\times H$ for some measurable $B\subset\Tht$. Since the prior
is of product form, $\Pi=\Pi_\Tht\times\Pi_H$, the
marginal posterior for the parameter ${\theta}\in\Tht$ depends on the
nuisance factor only through the integrated likelihood ratio,
%
%e2.5 #&#
%
\begin{equation}
\label{eqdefSn}
S_n\dvtx\Tht\rightarrow\RR\dvtx {\theta}\mapsto\int_H\prod_{i=1}^n
\frac{p_{{\theta},\eta}}{p_{{\theta}_0,\eta_0}}(X_i) \,d\Pi_H(\eta),
\end{equation}
where we have introduced factors $p_{{\theta}_0,\eta_0}(X_i)$ in the
denominator
for later convenience; see~(\ref{eqDefPD}). [The localized version of
(\ref{eqdefSn}) is denoted $h\mapsto s_n(h)$; see~(\ref{eqdefsn}).]
The map $S_n$ is to be viewed in a role similar to
that of the \textit{profile likelihood} in semiparametric maximum-likelihood
methods (see, e.g., Severini and Wong~\cite{Severini92} and
Murphy and van der Vaart~\cite{Murphy00}), in the sense
that~$S_n$ embodies the intermediate stage between nonparametric and
semiparametric steps of the estimation procedure.

We impose smoothness through a form of Le Cam's local
asymptotic normality: let $P\in\scrP$ be given, and let
$t\mapsto P_t$ be a one-dimensional submodel of $\scrP$ such that
$P_{t=0}=P$. Specializing to i.i.d. observations, we say that the model
is \textit{stochastically LAN} at $P\in\scrP$ along the direction
$t\mapsto P_t$, if there exists an $L_2(P)$-function $g_P$ with $Pg_P=0$
such that for all random sequences~$(h_n)$ bounded in $P$-probability,
%
%e2.6 #&#
%
\begin{equation}
\label{eqslan}\quad
\log\prod_{i=1}^n\frac{p_{n^{-1/2}h_n}}{p}(X_i) = \frac{1}{\sqrt{n}}
\sum_{i=1}^n h_n^Tg_P(X_i)-\ft12 h_n^T I_P h_n + o_P(1).
\end{equation}
Here $g_P$ is the score-function, and $I_P=P(g_P)^2$ is the Fisher
information of the submodel at $P$. Stochastic LAN is slightly
stronger than the usual LAN property~\cite{LeCam53,LeCam90}. In
examples, the proof of the ordinary LAN property often extends
to stochastic LAN without significant difficulties.

\begin{figure}

\includegraphics{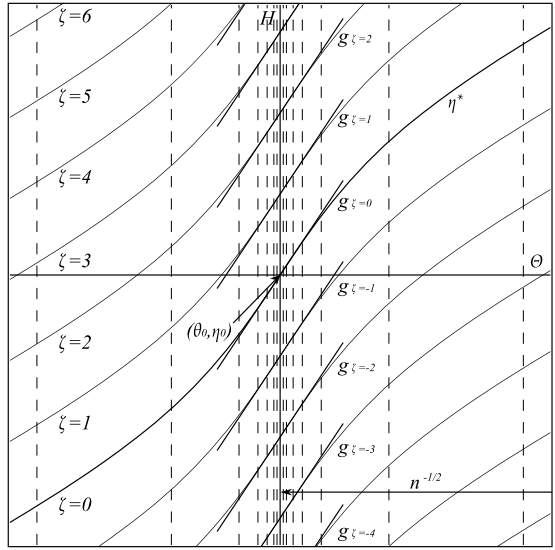}

\caption{A neighborhood of $({\theta}_0,\eta_0)$. Curved lines
represent sets $\{({\theta},\zeta)\dvtx{\theta}\in U_0\}$ for fixed~$\zeta$. The
curve through $\zeta=0$ parametrizes the least-favorable
submodel. Vertical dashed lines delimit regions such that
$\|{\theta}-{\theta}_0\|\leq n^{-1/2}$. Also indicated
are directions along which the likelihood is expanded, with score
functions $g_\zeta$.}
\label{figrecoord}
\end{figure}

Although formally only a convenience, the presentation benefits from
an \textit{adaptive} reparametrization (see Section 2.4 of
Bickel et al.~\cite{Bickel98}): based on the least-favorable submodel
$\eta^*$, we define, for all ${\theta}\in U_0$, $\eta\in H$,
%
%e2.7 #&#
%
\begin{equation}
\label{eqrepara}
({\theta},\eta({\theta},\zeta)) = \bigl({\theta},\eta^*({\theta
})+\zeta\bigr),\qquad
({\theta},\zeta({\theta},\eta)) = \bigl({\theta},\eta-\eta
^*({\theta})\bigr),
\end{equation}
and we introduce the notation $Q_{{\theta},\zeta}=P_{{\theta},\eta
({\theta},\zeta)}$.
With $\zeta=0$, ${\theta}\mapsto Q_{{\theta},0}$ describes the
least-favorable
submodel $\scrP^*$ and with a nonzero value of $\zeta$,
${\theta}\mapsto Q_{{\theta},\zeta}$ describes a version thereof, translated
over a nuisance direction (see Figure~\ref{figrecoord}).
Expressed in terms of the metric $r_H(\zeta_1,\zeta_2)=H(Q_{{\theta
}_0,\zeta_1},
Q_{{\theta}_0,\zeta_2})$, the sets $D({\theta},\rho)$ are mapped to
open balls
$B(\rho)=\{\zeta\in H\dvtx r_H(\zeta,0)<\rho\}$ centered at the
origin $\zeta=0$,
\[
\{ P_{{\theta},\eta} \dvtx {\theta}\in U_0,\eta\in D({\theta},\rho) \}
= \{ Q_{{\theta},\zeta} \dvtx {\theta}\in U_0,\zeta\in B(\rho)\}.
\]
In the formulation of Theorem~\ref{thmsbvmone}, we make use of a
domination condition based on the quantities
\[
U_n(\rho,h) = \sup_{\zeta\in B(\rho)} Q_{{\theta}_0,\zeta
}^n\Biggl(
\prod_{i=1}^n \frac{q_{{\theta}_n(h),\zeta}}
{q_{{\theta}_0,\zeta}}(X_i)\Biggr)
\]
for all $\rho>0$ and $h\in\RR^k$. Below, it is required
that there exists a sequence~$(\rho_n)$ with $\rho_n\downarrow0$,
$n\rho_n^2\rightarrow\infty$, such that, for every \textit{bounded},
stochastic sequence~$(h_n)$, $U(\rho_n,h_n)=O(1)$ (where the
expectation concerns the stochastic dependence of $h_n$ as well;
see \textit{Notation and conventions}). For a single, fixed
$\zeta$, the requirement says that the likelihood ratio remains
integrable when we replace ${\theta}_n(h_n)$ by the maximum-likelihood
estimator $\hat{{\theta}}_n(X_1,\ldots,X_n)$. Lemma~\ref{lemUdom}
demonstrates that ordinary differentiability of the likelihood-ratio
with respect to $h$, combined with a uniform upper bound on certain
Fisher information coefficients, suffices to satisfy $U(\rho_n,h_n)=O(1)$
for all bounded, stochastic $(h_n)$ and every $\rho_n\downarrow0$.

\begin{figure}

\includegraphics{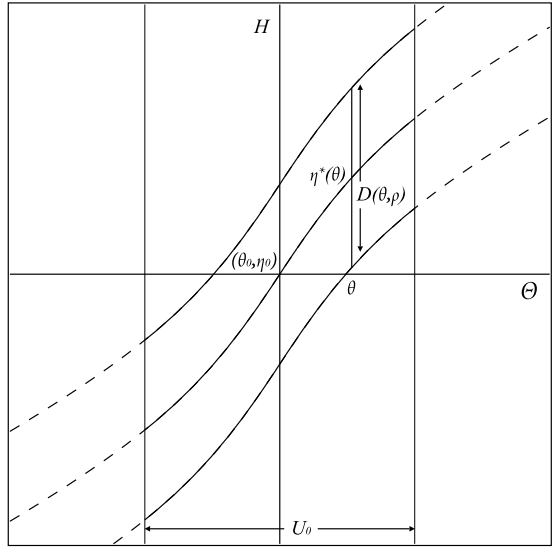}

\caption{A neighborhood of $({\theta}_0,\eta_0)$.
Shown are the least-favorable curve $\{({\theta},\eta^*({\theta
}))\dvtx\allowbreak{\theta}\in
U_0\}$
and (for fixed ${\theta}$ and $\rho>0$) the neighborhood $D({\theta
},\rho)$ of
$\eta^*({\theta})$. The sets $D({\theta},\rho)$ are expected to capture
(${\theta}$-conditional) posterior mass one asymptotically, for all
$\rho>0$ and ${\theta}\in U_0$.}
\label{fignbd-D}
\end{figure}

The second step of the proof can now be summarized as follows:
assuming stochastic LAN of the model, contraction of the nuisance
posterior as in Figure~\ref{fignbd-D} and said domination condition
%
%f1 #&#
%
are enough to turn LAN expansions for the integrand in~(\ref{eqdefSn})
into a single LAN expansion for $S_n$. The latter is determined by the
efficient score, because the locus of posterior concentration,~$\scrP^*$,
is a least-favorable submodel (see Theorem~\ref{thmilanone}).

%
%f2 #&#
%

The third step is based on two observations: first, in a semiparametric
problem, the integrals $S_n$ appear in the expression for the marginal
posterior in exactly the same way as parametric likelihood ratios
appear in the posterior for parametric problems. Second, the parametric
Bernstein--von Mises proof depends on likelihood ratios \textit{only}
through the LAN property. As a consequence, local asymptotic normality
for $S_n$ offers the possibility to apply Le Cam's proof of posterior
asymptotic normality in semiparametric context. If, in addition, we
impose contraction at parametric rate for the marginal posterior, the
LAN expansion of $S_n$ leads to the conclusion that the marginal
posterior satisfies the Bernstein--von Mises
assertion~(\ref{eqassertBvM}); see Theorem~\ref{thmpan}.

%s2.2 #&#
\subsection{Main theorem}
\label{submainthm}

Before we state the main result of this paper, general conditions
imposed on models and priors are formulated:
\begin{longlist}[(ii)]
\item[(i)] \textit{Model assumptions.}
Throughout the remainder of this article, $\scrP$ is assumed to be
well specified and dominated by a $\sigma$-finite measure on the
sample space and parametrized identifiably on $\Tht\times H$, with
$\Tht\subset\RR^k$ open and $H$ a subset of a metric vector-space
with metric $d_H$. Smoothness of the model is required but mentioned
explicitly throughout. We also assume that there exists an open
neighborhood $U_0\subset\Tht$ of ${\theta}_0$ on which a least-favorable
submodel $\eta^*\dvtx U_0\rightarrow H$ is defined.
%such that $\tht\mapsto P_{\tht,\eta^*(\tht)+\zeta}$ is stochastically
%LAN in the $\tht$-direction, for all $\zeta$ in an open neighborhood
%of $\zeta=0$. Furthermore, we assume that the efficient Fisher
%information at $(\tht_0,\eta_0)$ is nonsingular.
%
\item[(ii)] \textit{Prior assumptions.}
With regard to the prior $\Pi$ we follow the product structure of
the parametrization of $\scrP$, by endowing the parameterspace
$\Tht\times H$ with a product-prior $\Pi_\Tht\times\Pi_H$ defined
on a $\sigma$-field that includes the Borel $\sigma$-field generated
by the product-topology. Also, it is assumed that the prior $\Pi_\Tht$
is thick at~${\theta}_0$.
\end{longlist}
With the above general considerations for model and prior in mind, we
formulate the main result of this paper.
%
%th2.1 #&#
%
\begin{Thm}[(Semiparametric Bernstein--von Mises)]
\label{thmsbvmone}
Let $X_1,X_2,\ldots$ be distributed i.i.d.-$P_0$, with $P_0\in
\scrP$, and
let $\Pi_\Tht$ be thick at ${\theta}_0$. Suppose that for large
enough $n$,
the map $h\mapsto s_n(h)$ is continuous $P_0^n$-almost-surely. Also
assume that ${\theta}\mapsto Q_{{\theta},\zeta}$ is stochastically
LAN in the
${\theta}$-direction, for all $\zeta$ in an $r_H$-neighborhood of
$\zeta=0$ and that the efficient Fisher information $\effFI_{{\theta
}_0.\eta_0}$
is nonsingular. Furthermore, assume that there exists a sequence
$(\rho_n)$ with $\rho_n\downarrow0$, $n\rho_n^2\rightarrow\infty$
such that:
\begin{longlist}[(ii)]
\item[(i)] For all $M>0$, there exists a $K>0$ such that, for large
enough $n$,
\[
\Pi_H( K_n(\rho_n,M) ) \geq e^{-Kn\rho_n^2}.
\]
\item[(ii)] For all $n$ large enough, the Hellinger metric entropy
satisfies
\[
N(\rho_n,H,d_H)\leq e^{n\rho_n^2}
\]
and, for every bounded, stochastic $(h_n)$.
\item[(iii)] The model satisfies the domination condition,
%
%e2.8 #&#
%
\begin{equation}
\label{eqnormdom}
U_n(\rho_n,h_n) = O(1).
\end{equation}
\item[(iv)] For all $L>0$, Hellinger distances satisfy the uniform bound,
\[
\sup_{\{\eta\in H\dvtx d_H(\eta,\eta_0)\geq L\rho_n\}}
\frac{H(P_{{\theta}_n(h_n),\eta},P_{{\theta}_0,\eta
})}{H(P_{{\theta}_0,\eta},P_0)}=o(1).
\]
Finally, suppose that
\item[(v)] for every $(M_n)$, $M_n\rightarrow\infty$, the posterior
satisfies
\[
\Pi_n( \|h\|\leq M_n\bigmv X_1,\ldots,X_n )\convprob
{P_0} 1.
\]
\end{longlist}
Then the sequence of marginal posteriors for ${\theta}$ converges in total
variation to a normal distribution,
%
%e2.9 #&#
%
\begin{equation}
\label{eqConvTV}
\sup_{A}\bigl|
\Pi_n( h\in A\bigmv X_1,\ldots,X_n )
- N_{\effDelta_n,\effFI_{{\theta}_0,\eta_0}^{-1}}(A)
\bigr| \convprob{P_0} 0,
\end{equation}
centered on $\effDelta_n$ with covariance matrix $\effFI_{{\theta
}_0,\eta_0}^{-1}$.
\end{Thm}
\begin{pf}
The assertion follows from combination of Theorem~\ref{thmpertroc},
Corollary~\ref{corconspert}, Theorems~\ref{thmilanone} and
\ref{thmpan}.
\end{pf}

Let us briefly discuss some aspects of the conditions of
Theorem~\ref{thmsbvmone}. First, consider the required existence of
a least-favorable submodel in $\scrP$. In many semiparametric problems,
the efficient score function is \textit{not} a proper score in the sense
that it corresponds to a smooth submodel; instead, the efficient score
lies in the $L_2$-closure of the set of all proper scores. So there
exist sequences of so-called \textit{approximately least-favorable}
submodels whose scores converge to the efficient score in
$L_2$~\cite{vdVaart98}. Using such approximations of~$\scrP^*$, our
proof will entail extra conditions, but there is no reason to
expect problems of an overly restrictive nature. It may therefore be
hoped that the result remains largely unchanged
if we turn~(\ref{eqrepara}) into a sequence of reparametrizations
based on suitably chosen approximately least-favorable submodels.

Second, consider the rate $(\rho_n)$, which must be slow enough to
satisfy condition (iv) and is fixed at (or above) the minimax
Hellinger rate for estimation of the nuisance with known ${\theta}_0$
by condition (ii), while satisfying (i) and (iii)
as well. Conditions (i) and (ii) also arise
when considering Hellinger rates for nonparametric posterior convergence
and the methods of Ghosal et al.~\cite{Ghosal00} can be applied
in the present context with minor modifications. In addition,
Lemma~\ref{lemUdom} shows that in a wide class of semiparametric models,
condition (iii) is satisfied for \textit{any}
rate sequence $(\rho_n)$. Typically, the numerator in condition
(iv) is of order $O(n^{-1/2})$, so that condition (iv)
holds true for any $\rho_n$ such that $n\rho_n^2\rightarrow\infty$.
The above enables a rate-free version of the semiparametric
Bernstein--von Mises theorem (Corollary~\ref{corsimplesbvm}), in
which conditions (i) and (ii) above are weakened to become
comparable to those of Schwartz~\cite{Schwartz65}\vadjust{\goodbreak} for
nonparametric posterior consistency. Applicability of
Corollary~\ref{corsimplesbvm} is demonstrated in
Section~\ref{secplr}, where the linear coefficient in the
partial linear regression model is estimated.

Third, consider condition (v) of Theorem~\ref{thmsbvmone}:
though it is necessary [as it follows from~(\ref{eqConvTV})], it
is hard to formulate straightforward sufficient conditions to
satisfy (v) in generality. Moreover, condition (v)
involves the nuisance prior and, as such, imposes another condition
on $\Pi_H$ besides (i). To lessen its influence on~$\Pi_H$,
constructions in Section~\ref{secmarg} either work for all nuisance
priors (see Lem\-ma~\ref{lemlehmann}) or require only consistency
of the nuisance posterior (see Theorem~\ref{thmrocBayes}). The
latter is based on the limiting behavior of posteriors in misspecified
parametric models~\cite{Kleijn03,Kleijn07} and allows for the
tentative but general observation that a bias [cf.~(\ref{eqnearstraight})]
may ruin $n^{-1/2}$-consistency of the marginal posterior, especially
if the rate $(\rho_n)$ is sub-optimal. In the example of
Section~\ref{secplr}, the ``hard work'' stems from condition
(v) of
Theorem~\ref{thmsbvmone}: $\alpha>1/2$ H\"older smoothness and
boundedness of the family of regression functions in
Corollary~\ref{corsmoothplr} are imposed in order to satisfy
this condition. Since conditions (i) and (ii) appear
quite reasonable and conditions (iii) and (iv) are
satisfied relatively easily, condition (v) should be viewed
as the most complicated in an essential way.

To conclude, consistency under perturbation (with appropriate rate) is
one of the sufficient conditions, but it is by no means clear in how
far it should also hold with necessity. One expects that in some
situations where consistency under perturbation fails to hold fully,
integral local asymptotic normality (see Section~\ref{secilan}) is
still satisfied in a weaker form. In particular, it is possible that
(\ref{eqilan}) holds with a less-than-efficient score and Fisher
information, a result that would have an interpretation analogous
to suboptimality in H\'ajek's convolution theorem. What happens in
cases where integral LAN fails more comprehensively is both
interesting and completely mysterious from the point of view taken
in this article.

%%%%%%%%%%%%%%%%%%%%%%%%%%%%%%%%%%%%%%%%%%%%%%%%%%%%%%%%%%%%%%%%%%%%%%%%%%%%%%%%%%

%s3 #&#
\section{Posterior convergence under perturbation}
\label{secpert}

In this section, we consider contraction of the posterior around
least-favorable submodels. We express this form of posterior
convergence by showing that (under suitable conditions) the
conditional posterior for the nuisance parameter contracts around
the least-favorable submodel, conditioned on a sequence
${\theta}_n(h_n)$ for the parameter of interest with $h_n=O_{P_o}(1)$.
We view the sequence of models~$\scrP_{{\theta}_n(h_n)}$ as a random
perturbation of the model~$\scrP_{{\theta}_0}$ and generalize Ghosal
et al.~\cite{Ghosal00} to describe posterior contraction.
Ultimately, random perturbation of ${\theta}$ represents the
``appropriate form of uniformity'' referred to just after
definition~(\ref{eqKsets}). Given a rate sequence $(\rho_n)$,
$\rho_n\downarrow0$, we say that the conditioned nuisance
posterior is \textit{consistent under $n^{-1/2}$-perturbation at
rate} $\rho_n$, if
%
%e3.1 #&#
%
\begin{equation}
\label{eqpertroc}
\Pi_n\bigl( D^c({\theta},\rho_n) \bigmv
{\theta}={\theta}_0+n^{-1/2}h_n ; X_1,\ldots,X_n \bigr)\convprob
{P_0} 0
\end{equation}
for all bounded, stochastic sequences $(h_n)$.\vadjust{\goodbreak}
%
%th3.1 #&#
%
\begin{Thm}[(Posterior rate of convergence under
perturbation)]
\label{thmpertroc}
Assume that there exists a sequence $(\rho_n)$ with $\rho_n\downarrow0$,
$n\rho_n^2\rightarrow\infty$ such that for all $M>0$ and every bounded,
stochastic $(h_n)$:
\begin{longlist}[(iii)]
\item[(i)] There exists a constant $K>0$ such that for large enough $n$,
%
%e3.2 #&#
%
\begin{equation}
\label{eqsuffprior}
\Pi_H( K_n(\rho_n,M) ) \geq e^{-Kn\rho_n^2}.
\end{equation}
\item[(ii)] For $L>0$ large enough, there exist $(\phi_n)$ such
that for large enough~$n$,
%
%e3.3 #&#
%
\begin{equation}
\label{eqnuistest}
P_0^n\phi_n\rightarrow0,\qquad
\sup_{\eta\in D^c({\theta}_0,L\rho_n)}
P_{{\theta}_n(h_n),\eta}^n(1-\phi_n)\leq e^{-L^2n\rho_n^2/4}.
\end{equation}
\item[(iii)] The least-favorable submodel
satisfies $d_H(\eta^*({\theta}_n(h_n)),\eta_0)=o(\rho_n)$.
\end{longlist}
Then, for every bounded, stochastic $(h_n)$ there exists an $L>0$
such that the conditional nuisance posterior converges as
%
%e3.4 #&#
%
\begin{equation}
\label{eqrocassert}
\Pi\bigl( D^c({\theta},L\rho_n)\bigmv
{\theta}={\theta}_0+n^{-1/2}h_n; X_1,\ldots,X_n \bigr)= o_{P_0}(1)
\end{equation}
under $n^{-1/2}$-perturbation.
\end{Thm}
\begin{pf} Let $(h_n)$ be a stochastic sequence bounded by $M$,
and let $0<C<1$ be given. Let $K$ and $(\rho_n)$ be as in
conditions (i) and (ii). Choose $L>4\sqrt{1+K+C}$ and
large enough to satisfy condition (ii) for some $(\phi_n)$. By
Lem\-ma~\ref{lemrocdenom}, the events
\[
A_n = \Biggl\{
\int_H\prod_{i=1}^n\frac{p_{{\theta}_n(h_n),\eta}}{p_{{\theta
}_0,\eta_0}}(X_i)
\,d\Pi_H(\eta)
\geq e^{-(1+C)n\rho_n^2} \Pi_H(K_n(\rho_n,M)) \Biggr\}
\]
satisfy $P_0^n(A_n^c)\rightarrow0$. Using also
the first limit in~(\ref{eqnuistest}), we then derive
\begin{eqnarray*}
&&
P_0^n\Pi\bigl( D^c({\theta},L\rho_n)\bigmv {\theta}={\theta}_n(h_n);
X_1,\ldots,X_n \bigr)\\
&&\qquad\leq
P_0^n\Pi\bigl( D^c({\theta},L\rho_n)\bigmv {\theta}={\theta}_n(h_n);
X_1,\ldots,X_n \bigr) 1_{A_n}
(1-\phi_n) + o(1)
\end{eqnarray*}
[even with random $(h_n)$, the posterior
$\Pi( \cdot|{\theta}={\theta}_n(h_n); X_1,\ldots,X_n
)\leq1$,
by definition~(\ref{eqposterior})]. The first term on the r.h.s.
can be bounded further by the definition of the events $A_n$,
\begin{eqnarray*}
&&P_0^n\Pi\bigl( D^c({\theta},L\rho_n)\bigmv {\theta}={\theta}_n;
X_1,\ldots,X_n \bigr) 1_{A_n}
(1-\phi_n)\\
&&\qquad\leq
\frac{e^{(1+C)n\rho_n^2}}{\Pi_H(K_n(\rho_n,M))}
P_0^n\Biggl( \int_{D^c({\theta}_n(h_n),L\rho_n)}\prod_{i=1}^n
\frac{p_{{\theta}_n(h_n),\eta}}{p_{{\theta}_0,\eta_0}}(X_i)
(1-\phi_n)
\,d\Pi_H\Biggr).
\end{eqnarray*}
Due to condition (iii) it follows that
%
%e3.5 #&#
%
\begin{equation}
\label{eqDs}
D\biggl({\theta}_0,\ft12 L \rho_n\biggr)\subset
\bigcap_{n\geq1} D({\theta}_n(h_n),L\rho_n)
\end{equation}
for large enough $n$. Therefore,
%
%e3.6 #&#
%
\begin{eqnarray}
\label{eqonedomain}
&&
P_0^n \int_{D^c({\theta}_n(h_n),L\rho_n)}\prod_{i=1}^n
\frac{p_{{\theta}_n(h_n),\eta}}{p_{{\theta}_0,\eta_0}}(X_i)
(1-\phi_n)
\,d\Pi_H(\eta)\nonumber\\[-8pt]\\[-8pt]
&&\qquad\leq\int_{D^c({\theta}_0,L\rho_n/2)}
P_{{\theta}_n(h_n),\eta}^n(1-\phi_n) \,d\Pi_H(\eta).\nonumber
\end{eqnarray}
Upon substitution of~(\ref{eqonedomain}) and with the use of the
second bound in~(\ref{eqnuistest}) and~(\ref{eqsuffprior}),
the choice we made earlier for $L$ proves the assertion.
\end{pf}

We conclude from the above that besides sufficiency of prior mass, the
crucial condition for consistency under perturbation is the existence
of a~test sequence $(\phi_n)$ satisfying~(\ref{eqnuistest}). To find
sufficient conditions, we follow a~construction of tests
based on the Hellinger geometry\vspace*{1pt} of the model, generalizing the approach
of Birg\'e~\cite{Birge83,Birge84} and Le Cam~\cite{LeCam86} to
$n^{-1/2}$-perturbed context. It is easiest to illustrate their
approach by considering the problem of testing/estimating
$\eta$ when ${\theta}_0$ is known: we cover the nuisance model
$\{P_{{\theta}_0,\eta}\dvtx\eta\in H\}$ by a minimal collection of Hellinger
balls $B$ of radii $(\rho_n)$, each of which is convex and hence
testable against $P_0$ with power bounded by $\exp(-\ft14 n H^2(P_0,B))$,
based on the minimax theorem~\cite{LeCam86}. The tests for the covering
Hellinger balls are combined into a single test for the nonconvex
alternative $\{P\dvtx H(P,P_0)\geq\rho_n\}$ against $P_0$. The order of
the cover controls the power of the combined test. Therefore the
construction requires an upper bound to Hellinger metric entropy
numbers~\cite{vdVaart96}
%
%e3.7 #&#
%
\begin{equation}
\label{eqminimaxrate}
N(\rho_n,\scrP_{{\theta}_0},H) \leq e^{n\rho_n^2},
\end{equation}
which is interpreted as indicative of the nuisance model's complexity
in the sense that the lower bound to the collection of rates $(\rho_n)$
solving~(\ref{eqminimaxrate}) is the Hellinger minimax rate for
estimation of $\eta_0$. In the $n^{-1/2}$-perturbed
problem, the alternative does not just consist of the complement of a
Hellinger-ball in the nuisance factor $H$, but also has an extent
in the ${\theta}$-direction shrinking at rate $n^{-1/2}$. Condition
(\ref{eqHcone}) below guarantees that Hellinger covers of~$H$ like the above
are large enough to accommodate the ${\theta}$-extent of the alternative,
the implication being that the test sequence one constructs
for the nuisance in case ${\theta}_0$ is known, can also be used when
${\theta}_0$ is known only up to $n^{-1/2}$-perturbation. Therefore, the
entropy bound in Lemma~\ref{lemtestpert} is~(\ref{eqminimaxrate}).
Geometrically,~(\ref{eqHcone}) requires that $n^{-1/2}$-perturbed versions
of the nuisance model are contained in a narrowing sequence of metric
cones based at $P_0$. In differentiable models,\vspace*{1pt} the Hellinger distance
$H(P_{{\theta}_n(h_n),\eta},P_{{\theta}_0,\eta})$ is typically of
order $O(n^{-1/2})$
for all $\eta\in H$. So if, in addition, $n\rho_n^2\rightarrow\infty$,
limit~(\ref{eqHcone}) is expected to hold pointwise in $\eta$. Then
only the uniform character of~(\ref{eqHcone}) truly forms a condition.
%
%le3.2 #&#
%
\begin{lemma}[(Testing under perturbation)]
\label{lemtestpert}
If $(\rho_n)$ satisfies $\rho_n\downarrow0$,]\break
$n\rho_n^2\rightarrow\infty$ and the following requirements are met:
\begin{longlist}[(ii)]
\item[(i)] For all $n$ large enough, $N(\rho_n,H,d_H)
\leq e^{n\rho_n^2}$.
\item[(ii)] For all $L>0$ and all bounded, stochastic $(h_n)$,
%
%e3.8 #&#
%
\begin{equation}
\label{eqHcone}
\sup_{\{\eta\in H\dvtx d_H(\eta,\eta_0)\geq L\rho_n\}}
\frac{H(P_{{\theta}_n(h_n),\eta},P_{{\theta}_0,\eta
})}{H(P_{{\theta}_0,\eta},P_0)}=o(1).
\end{equation}
\end{longlist}
Then for all $L\geq4$, there exists a test sequence $(\phi_n)$ such
that for all bounded, stochastic $(h_n)$,
%
%e3.9 #&#
%
\begin{equation}
\label{eqDtests}
P_0^n\phi_n\rightarrow0,\qquad
\sup_{\eta\in D^c({\theta}_0,L\rho_n)}
P_{{\theta}_n(h_n),\eta}^n(1-\phi_n)\leq e^{-L^2n\rho_n^2/4}
\end{equation}
for large enough $n$.
\end{lemma}
\begin{pf} Let $(\rho_n)$ be such that (i) and (ii)
are satisfied. Let $(h_n)$ and $L\geq4$ be given. For all $j\geq1$,
define $H_{j,n}=\{\eta\in H\dvtx jL\rho_n\leq d_H(\eta_0,\eta)\leq
(j+1)L\rho
_n\}$
and $\scrP_{j,n}=\{P_{{\theta}_0,\eta}\dvtx\eta\in H_{j,n}\}$. Cover
$\scrP_{j,n}$
with Hellinger balls $B_{i,j,n}(\ft14jL\rho_n)$, where
\[
B_{i,j,n}(r)=\{ P\dvtx H(P_{i,j,n},P)\leq r \}
\]
and $P_{i.j.n}\in\scrP_{j,n}$, that is, there exists an $\eta
_{i,j,n}\in
H_{j,n}$ such that $P_{i,j,n}=P_{{\theta}_0,\eta_{i,j,n}}$. Denote
$H_{i,j,n}=\{\eta\in H_{j,n}\dvtx P_{{\theta}_0,\eta}\in B_{i,j,n}
(\ft14jL\rho_n)\}$. By assumption, the minimal number of
such balls needed to cover $\scrP_{i,j}$ is finite; we denote the
corresponding covering number by $N_{j,n}$, that is, $1\leq i\leq N_{j,n}$.

Let $\eta\in H_{j,n}$ be given. There exists an $i$ ($1\leq i\leq
N_{j,n}$) such that $d_H(\eta,\break\eta_{i,j,n})\leq\ft14 jL\rho_n$.
Then,\vspace*{1pt} by the triangle inequality, the definition of
$H_{j,n}$ and assumption~(\ref{eqHcone}),
%
%e3.10 #&#
%
\begin{eqnarray}
\label{eqHgeometry}
&&
H\bigl(P_{{\theta}_n(h_n),\eta},P_{{\theta}_0,\eta_{i,j,n}}\bigr)\nonumber\hspace*{-30pt}\\
&&\qquad\leq H\bigl(P_{{\theta}_n(h_n),\eta},P_{{\theta}_0,\eta}\bigr)
+ H(P_{{\theta}_0,\eta},P_{{\theta}_0,\eta_{i,j,n}})\nonumber\hspace*{-30pt}\\
&&\qquad\leq \frac{H(P_{{\theta}_n(h_n),\eta},P_{{\theta}_0,\eta})}
{H(P_{{\theta}_0,\eta},P_0)} H(P_{{\theta}_0,\eta},P_0)+ \ft14 jL\rho_n\hspace*{-30pt}\\
&&\qquad\leq\biggl(\sup_{\{\eta\in H\dvtx d_H(\eta,\eta_0)\geq L\rho_n\}}
\frac{H(P_{{\theta}_n(h_n),\eta},P_{{\theta}_0,\eta})}{H(P_{{\theta}_0,\eta
},P_0)}\biggr)(j+1)L\rho_n + \ft14 jL\rho_n\nonumber\hspace*{-30pt}\\
&&\qquad\leq\ft12 jL\rho_n\nonumber\hspace*{-30pt}
\end{eqnarray}
for large enough $n$. We conclude that there exists an
$N\geq1$ such that for all $n\geq N$,
$j\geq1$, $1\leq i\leq N_{j,n}$, $\eta\in H_{i,j,n}$,
$P_{{\theta}_n(h_n),\eta}\in B_{i,j,n}(\ft12jL\rho_n)$. Moreover, Hellinger
balls are convex and for all $P\in B_{i,j,n}(\ft12jL\rho_n)$,
$H(P,P_0)\geq\ft12jL\rho_n$. As a consequence of the minimax theorem
(see Le Cam~\cite{LeCam86}, Birg\'e
\cite{Birge83,Birge84}), there exists a test sequence
$(\phi_{i,j,n})_{n\geq1}$ such that
\[
P_0^n\phi_{i,j,n} \vee \sup_{P}
P^n(1-\phi_{i,j,n})
\leq e^{-nH^2(B_{i,j,n}(jL\rho_n/2),P_0)}
\leq e^{-nj^2L^2\rho_n^2/4},
\]
where the\vspace*{1pt} supremum runs over all $P\in B_{i,j,n}(\frac12jL\rho_n)$.
Defining, for all $n\geq1$,
$\phi_n = \sup_{j\geq1}\max_{1\leq i\leq N_{j,n}} \phi_{i,j,n}$,
we find (for details,\vspace*{1pt} see the proof of Theorem~3.10 in
\cite{Kleijn03}) that
%
%e3.11 #&#
%
\begin{equation}
\label{eqtestintermediate}
P_0^n\phi_n \leq\sum_{j\geq1} N_{j,n} e^{-L^2j^2n\rho_n^2/4},\qquad
P^n(1-\phi_n) \leq e^{-L^2n\rho_n^2/4}
\end{equation}
for all $P=P_{{\theta}_n(h_n),\eta}$ and $\eta\in D^c({\theta
}_0,L\rho_n)$.
Since $L\geq4$, we have for all $j\geq1$,
%
%e3.12 #&#
%
\begin{eqnarray}
\label{eqestN}
N_{j,n} &=& N\bigl(\tfrac14Lj\rho_n,\scrP_{j,n},H\bigr)
\leq N\bigl(\tfrac14Lj\rho_n,\scrP,H\bigr)\nonumber\\[-8pt]\\[-8pt]
&\leq& N(\rho_n,\scrP,H) \leq e^{n\rho_n^2}\nonumber
\end{eqnarray}
by assumption~(\ref{eqminimaxrate}). Upon substitution of
(\ref{eqestN}) into~(\ref{eqtestintermediate}), we obtain the
following bounds:
\[
P_0^n\phi_n\leq\frac{e^{(1-L^2/4)n\rho_n^2}}
{1-e^{-L^2n\rho_n^2/4}},\qquad
\sup_{\eta\in D^c({\theta}_0,L\rho_n)}
P_{{\theta}_n(h_n),\eta}^n(1-\phi_n)\leq e^{-L^2n\rho_n^2/4}
\]
for large enough $n$, which implies assertion~(\ref{eqDtests}).
\end{pf}

%For some models, the sequence of bounds~(\ref{eqestN}) is too
%coarse. Problems arise already for finite-dimensional
%parameter spaces if they are unbounded: while the \lhs\ of
%(\ref{eqestN}) is finite, subsequent bounds are infinite because
%totally-boundedness is lost. In such cases, we would forgo
%estimations~(\ref{eqestN}) and control $N_{j,n}$ more directly.
%
%Possible generalization of theorem~\ref{thmpertroc}
%relates to the size of the perturbation. Since we apply
%theorem~\ref{thmpertroc} only in differentiable situations, we
%specialize the proof here to perturbations of size $n^{-1/2}$ and
%rely on differentiability to achieve inclusion~(\ref{eqDs}).
%However, if we can achieve~(\ref{eqDs}) in another way, the argument
%based on~(\ref{eqHgeometry}) shows that the construction given
%above can be generalized to perturbations of any size $\tau_n$
%such that $\tau_n = o(\rho_n)$. This would enable study of
%consistency and rates of convergence under perturbations of
%larger than parametric order, which appears most appealing in
%situations where the full, nonparametric posterior is known to
%converge at rate $\tau_n$: in that case, the above would further
%specify posterior concentration to occur around $\eta^*$ at any
%rate $\rho_n$ above $\tau_n$. Such a generalization appears useful
%when the stochastic LAN expansion of the likelihood hinges on
%a rate different from $n^{-1/2}$ (for an example, see Kleijn and
%Knapik~\cite{Kleijn10a}).
%
In preparation of Corollary~\ref{corsimplesbvm}, we also provide
a version of Theorem~\ref{thmpertroc} that only asserts
consistency under $n^{-1/2}$-perturbation at \textit{some} rate
while relaxing bounds for prior mass and entropy. In the statement
of the corollary, we make use of the family of Kullback--Leibler
neighborhoods that would play a role for the posterior of the
nuisance if ${\theta}_0$ were known~\cite{Ghosal00}.
%
%e3.13 #&#
%
\begin{equation}
\label{eqsimpleKsets}\quad
K(\rho) = \biggl\{ \eta\in H \dvtx
-P_0\log\frac{p_{{\theta}_0,\eta}}{p_{{\theta}_0,\eta_0}}\leq
\rho^2,
P_0\biggl(\log\frac{p_{{\theta}_0,\eta}}{p_{{\theta}_0,\eta
_0}}\biggr)^2
\leq\rho^2 \biggr\}
\end{equation}
for all $\rho>0$. The proof below follows steps similar to those
in the proof of Corollary 2.1 in~\cite{Kleijn06}.
%
%co3.3 #&#
%
\begin{corollary}[(Posterior consistency under perturbation)]
\label{corconspert}
Assume that for all $\rho>0$, $N(\rho,H,d_H) < \infty$,
$\Pi_H( K(\rho)) > 0$ and:
\begin{longlist}[(ii)]
\item[(i)] For all $M>0$ there is an $L>0$ such that for all
$\rho>0$ and large enough $n$, $K(\rho) \subset K_n(L\rho,M)$.
\item[(ii)] For every bounded random sequence $(h_n)$,
$\sup_{\eta\in H} H(P_{{\theta}_n(h_n),\eta},P_{{\theta}_0,\eta
})$ and
$H(P_{{\theta}_0,\eta^*({\theta}_n(h_n))},P_{{\theta}_0,\eta_0})$
are of order
$O(n^{-1/2})$.
\end{longlist}
Then there exists a sequence
$(\rho_n)$, $\rho_n\downarrow0$, $n\rho_n^2\rightarrow\infty$, such
that the conditional nuisance posterior converges under
$n^{-1/2}$-perturbation at rate $(\rho_n)$.\vadjust{\goodbreak}
\end{corollary}
\begin{pf}
%Define functions $g_1$, $g_2$ and $g_n$,
% g_1(\rho)=\Pi_H( K(\rho) ),
% g_2(\rho)= N( \rho, \scrP_{\tht_0}, H ),
% g_n(\rho)=e^{-n\rho^2}(g_1(\rho)+\frac{1}{g_2(\rho)}).
%For large enough $n$, the functions $g_n$ are well defined and finite
%by the assumptions and $g_n(\rho)\rightarrow0$
%as $n\rightarrow\infty$, for every fixed $\rho>0$. Therefore, there
%exists a sequence $(\rho_n)$ such that $\rho_n\downarrow0$ and
%$n\rho_n^2\rightarrow\infty$, with $g_n(\rho_n)\rightarrow0$
%(\eg\ fix $n_1<n_2<\cdots$ large enough, such that $g_n(1/k)\leq1/k$
%for all $n\geq n_k$; next define $\rho_n=1/k$ for $n_k\leq n<
%n_{k+1}$).
%In particular, there exists an $N$ such that $g_n(\rho_n)\leq1$
%for all $n\geq N$. This implies that for all $n$ large enough,
%$g_1(\rho_n)\geq e^{-n\rho_n^2}$, so that~(\ref{eqsuffprior}) is
%satisfied, and $g_2(\rho_n)\leq e^{n\rho_n^2}$, so that
%(\ref{eqminimaxrate}) is satisfied.
We follow the proof of Corollary 2.1 in Kleijn and van der
Vaart~\cite{Kleijn06} and add that, under condition (ii),
(\ref{eqHcone}) and condition (iii) of Theorem~\ref{thmpertroc}
are satisfied. We conclude that there exists a test sequence
satisfying~(\ref{eqnuistest}). Then the assertion of
Theorem~\ref{thmpertroc} holds.
\end{pf}

The following lemma generalizes Lemma 8.1 in Ghosal et al.
\cite{Ghosal00} to the $n^{-1/2}$-perturbed setting.
%
%le3.4 #&#
%
\begin{lemma}
\label{lemrocdenom}
Let $(h_n)$ be stochastic and bounded by some $M>0$. Then
%
%e3.14 #&#
%
\begin{eqnarray}
\label{eqrocdenom}\qquad
&&
P_0^n\Biggl( \int_H\prod_{i=1}^n
\frac{p_{{\theta}_n(h_n),\eta}}{p_{{\theta}_0,\eta_0}}(X_i) \,d\Pi
_H(\eta)
< e^{-(1+C)n\rho^2} \Pi_H(K_n(\rho,M)) \Biggr)\nonumber\\[-10pt]\\[-10pt]
&&\qquad\leq\frac{1}{C^2n\rho^2}\nonumber
\end{eqnarray}
for all $C>0$, $\rho>0$ and $n\geq1$.\vspace*{-3pt}
\end{lemma}
\begin{pf} See the proof of Lemma 8.1 in Ghosal et al.
\cite{Ghosal00} (dominating the $h_n$-dependent log-likelihood ratio
immediately after the first application of Jensen's inequality).
%[OLD NOTATION!!]
%Let $(h_n)$, $C>0$, $\rho>0$ and $n\geq1$ be given.
%If $\Pi_H(K_{\rho,n})=0$, then~(\ref{eqrocdenom}) holds,
%so we assume $\Pi_H(K_{\rho,n})>0$ without loss of generality and
%denote the conditional prior by $\Pi_n(A)=\Pi_H( A | K_{\rho,n} )$
%for
%all measurable $A\subset H$. Then,
% \begin{split}
% P_0^n( \int_H\prod_{i=1}^n
% \frac{p_{\tht_n(h_n),\eta}}{p_{\tht_0,\eta_0}}(X_i)& d\Pi_H(\eta)
% < e^{-(1+C)n\rho^2} \Pi_H(K_{\rho,n}) )\\
% &\leq P_0^n( \int\prod_{i=1}^n
% \frac{p_{\tht_n(h_n),\eta}}{p_{\tht_0,\eta_0}}(X_i)
% d\Pi_n(\eta)< e^{-(1+C)n\rho^2} ).
% \end{split}
%By Jensen's inequality and~(\ref{eqKsets}),
% \log\int\prod_{i=1}^n
% \frac{p_{\tht_n(h_n),\eta}}{p_{\tht_0,\eta_0}}(X_i) d\Pi_n(\eta)
% \geq\sqrt{n}\int\GG_n\log\frac{p_{\tht_n(h_n),\eta}}{p_{\tht_0,
% d\Pi_n(\eta) - n\rho^2.
%so that,
% \begin{split}
% P_0^n( \int
% \prod_{i=1}^n \frac{p_{\tht_n(h_n),\eta}}{p_{\tht_0,\eta_0}}(X_i)
% & d\Pi_n(\eta)< e^{-(1+C)n\rho^2} )\\
% &\leq
% P_0^n( \int\GG_n\log\frac{p_{\tht_n(h_n),\eta}}{p_{\tht_0,
% d\Pi_n(\eta)
% < -\sqrt{n}C\rho^2 ).
% \end{split}
%By Chebyshev's inequality, Jensen's inequality, Fubini's theorem
%and the fact that for any sequence $(Z_n)$ in $L_2(P_0)$,
%$P_0^n(\GG_nZ_n)^2 ={\rm Var}_{P_0}Z_n\leq P_0^nZ_n^2$,
% \begin{split}
% P_0^n( \int
% \GG_n\log\frac{p_{\tht_n(h_n),\eta}}{p_{\tht_0,\eta_0}}& d\Pi_n(\eta)
% < -\sqrt{n}C\rho^2 )\\
% &\leq
% \frac{1}{nC^2\rho^4}\int P_0^n(
% \GG_n\log\frac{p_{\tht_n(h_n),\eta}}{p_{\tht_0,\eta_0}})^2
% d\Pi_n(\eta)
% \leq\frac{1}{nC^2\rho^2}.
% \end{split}
%where the last step follows from definition~(\ref{eqKsets}).
\end{pf}

%%%%%%%%%%%%%%%%%%%%%%%%%%%%%%%%%%%%%%%%%%%%%%%%%%%%%%%%%%%%%%%%%%%%%%%%%%%%%%%%%%

%s4 #&#
\section{Integrating local asymptotic normality}
\label{secilan}

The smoothness condition in the Le Cam's parametric Bernstein--von Mises
theorem is a LAN expansion of the likelihood, which is replaced in
semiparametric context by a stochastic LAN expansion of the integrated
likelihood~(\ref{eqdefSn}). In this section, we consider sufficient
conditions under which the localized integrated likelihood
%
%e4.1 #&#
%
\begin{equation}
\label{eqdefsn}
s_n(h) = \int_H \prod_{i=1}^n
\frac{p_{{\theta}_0+n^{-1/2}h,\eta}}{p_{{\theta}_0,\eta_0}}(X_i)
\,d\Pi_H(\eta)
\end{equation}
has the \textit{integral LAN} property; that is, $s_n$ allows an expansion
of the form
%
%e4.2 #&#
%
\begin{equation}
\label{eqilan}
\log\frac{s_n(h_n)}{s_n(0)} = \frac{1}{\sqrt{n}}\sum_{i=1}^\infty
h_n^T\effscore_{{\theta}_0,\eta_0}
-\ft12 h_n^T \effFI_{{\theta}_0,\eta_0} h_n + o_{P_0}(1)
\end{equation}
for every random sequence $(h_n)\subset\RR^k$ of order $O_{P_0}(1)$,
as required in Theorem~\ref{thmpan}. Theorem~\ref{thmilanone}
assumes that the model is stochastically LAN and requires consistency
under $n^{-1/2}$-perturbation for the nuisance posterior.
Consistency not only allows us to restrict sufficient
conditions to neighborhoods of $\eta_0$ in $H$, but also enables
lifting of the LAN expansion of the integrand in~(\ref{eqdefsn})
to an expansion of the integral $s_n$ itself; cf.~(\ref{eqilan}).
The posterior concentrates on the least-favorable submodel so that
only the least-favorable expansion at $\eta_0$ contributes to
(\ref{eqilan}) asymptotically. For this reason, the intergral LAN
expansion is determined by the efficient score function (and not some
other influence function). Ultimately, occurrence of the efficient
score lends the marginal posterior (and statistics based upon~it)
properties of frequentist semiparametric optimality.\vadjust{\goodbreak}

To derive Theorem~\ref{thmilanone}, we reparametrize the model; cf.
(\ref{eqrepara}). While yielding adaptivity, this reparametrization
also leads to ${\theta}$-dependence in the prior for~$\zeta$, a
technical issue that we tackle before addressing the main point of this
section. We show that the prior mass of the relevant neighborhoods
displays the appropriate type of stability, under a condition on local
behavior of Hellinger \mbox{distances} in the least-favorable model. For
smooth least-favorable submodels, typically
$d_H(\eta^*({\theta}_n(h_n)),\eta_0)=O(n^{-1/2})$ for all bounded,
stochastic $(h_n)$, which suffices.\vspace*{-3pt}
%
%le4.1 #&#
%
\begin{lemma}[(Prior stability)]
\label{lemtranslate}
Let $(h_n)$ be a bounded, stochastic sequence of perturbations, and
let $\Pi_H$ be any prior on $H$. Let $(\rho_n)$ be such that
$d_H(\eta^*({\theta}_n(h_n)),\eta_0) = o(\rho_n)$. Then the
prior mass of radius-$\rho_n$ neighborhoods of $\eta^*$ is
stable, that is,
%
%e4.3 #&#
%
\begin{equation}
\label{eqstab}
\Pi_H(D({\theta}_n(h_n),\rho_n))
=\Pi_H(D({\theta}_0,\rho_n)) + o(1).\vspace*{-3pt}
\end{equation}
\end{lemma}
\begin{pf}
Let $(h_n)$ and $(\rho_n)$ be such that
$d_H(\eta^*({\theta}_n(h_n)),\eta_0) = o(\rho_n)$.
Denote $D({\theta}_n(h_n),\rho_n)$ by $D_n$ and $D({\theta}_0,\rho
_n)$ by $C_n$
for all $n\geq1$. Since
\[
| \Pi_H(D_n) - \Pi_H(C_n) |
\leq\Pi_H\bigl((D_n\cup C_n)\setminus(D_n\cap C_n)\bigr),
\]
we consider the sequence of symmetric differences. Fix some
$0<\alpha<1$. Then for all $\eta\in D_n$ and all $n$ large enough,
$d_H(\eta,\eta_0)\leq d_H(\eta,\eta^*({\theta}_n(h_n))) +
d_H(\eta^*({\theta}_n(h_n)),\eta_0) \leq(1+\alpha)\rho_n$, so
that $D_n\cup C_n\subset D({\theta}_0,(1+\alpha)\rho_n)$.
Furthermore, for large enough $n$ and any $\eta\in
D({\theta}_0,(1-\alpha)\rho_n)$, $d_H(\eta,\eta^*({\theta}_n(h_n)))
\leq d_H(\eta,\eta_0) + d_H(\eta_0,\eta^*({\theta}_n(h_n)))\leq
\rho_n
+ d_H(\eta_0,\eta^*({\theta}_n(h_n)))-\alpha\rho_n<\rho_n$,
so that $D({\theta}_0,(1-\alpha)\rho_n)\subset D_n\cap C_n$.
Therefore,
\[
(D_n\cup C_n)\setminus(D_n\cap C_n) \subset
D\bigl({\theta}_0,(1+\alpha)\rho_n\bigr)\setminus D\bigl({\theta
}_0,(1-\alpha)\rho_n\bigr)
\rightarrow\varnothing,
\]
which implies~(\ref{eqstab}).\vspace*{-3pt}
\end{pf}

Once stability of the nuisance prior is established,
Theorem~\ref{thmilanone} hinges on stochastic local asymptotic
normality of the submodels $t\mapsto Q_{{\theta}_0+t,\zeta}$,
for all~$\zeta$ in an $r_H$-neighborhood of $\zeta=0$. We
assume there exists a $g_\zeta\in L_2(Q_{{\theta}_0,\zeta})$ such that
for every random $(h_n)$ bounded in $Q_{{\theta}_0,\zeta}$-probability,
%
%e4.4 #&#
%
\begin{equation}
\label{eqqlan}
\log\prod_{i=1}^n\frac{q_{{\theta}+n^{-1/2}h_n,\zeta}}{q_{{\theta
}_0,0}}(X_i)
= \frac{1}{\sqrt{n}}\sum_{i=1}^n h_n^Tg_\zeta(X_i)
-\ft12 h_n^TI_{\zeta}h_n+ R_n(h_n,\zeta),\hspace*{-35pt}
\end{equation}
where $I_{\zeta}=Q_{{\theta}_0,\zeta}g_\zeta g_\zeta^T$ and
$R_n(h_n,\zeta)=o_{Q_{{\theta}_0,\zeta}}(1)$. Equation~(\ref{eqqlan})
specifies the (minimal) tangent set (van der Vaart
\cite{vdVaart98}, Section 25.4) with\vspace*{1pt} respect to which
differentiability of the model is required. Note that
$g_0=\effscore_{{\theta}_0,\eta_0}$.\vspace*{-3pt}
%
%th4.2 #&#
%
\begin{Thm}[(Integral local asymptotic normality)]
\label{thmilanone}
Suppose that ${\theta}\mapsto Q_{{\theta},\zeta}$ is stochastically LAN
for all $\zeta$ in an $r_H$-neighborhood of
$\zeta=0$. Furthermore, assume that posterior consistency under
$n^{-1/2}$-perturbation obtains with a rate $(\rho_n)$ also valid
in~(\ref{eqnormdom}). Then the integral LAN-expansion~(\ref{eqilan})
holds.\vspace*{-3pt}\vadjust{\goodbreak}
\end{Thm}
\begin{pf} Throughout this proof $G_n(h,\zeta)=
\sqrt{n} h^T\PP_ng_\zeta-\frac12 h^T I_\zeta h$,
for all~$h$ and all $\zeta$. Furthermore, we abbreviate
${\theta}_n(h_n)$ to ${\theta}_n$ and omit explicit notation for
$(X_1,\ldots,X_n)$-dependence in several places.

Let $\delta,\ep>0$ be given, and let
${\theta}_n={\theta}_0+n^{-1/2}h_n$ with $(h_n)$ bounded in
$P_0$-probability. Then\vspace*{1pt}
there exists a constant $M>0$ such that $P_0^n(\|h_n\|>M)<\ft12\delta$
for all $n\geq1$. With $(h_n)$ bounded, the assumption of consistency
under $n^{-1/2}$-perturbation says that
\[
P_0^n\bigl(
\log\Pi\bigl( D({\theta},\rho_n)\bigmv
{\theta}={\theta}_n ; X_1,\ldots,X_n \bigr)\geq-\ep \bigr)
> 1-\tfrac12\delta
\]
for large enough $n$. This implies that the posterior's numerator
and denominator are related through
%
%e4.5 #&#
%
\begin{eqnarray}
\label{eqintD}
&&
P_0^n\Biggl(
\int_H \prod_{i=1}^n\frac{p_{{\theta}_n,\eta}}{p_{{\theta
}_0,\eta_0}}(X_i)
\,d\Pi_H(\eta)\nonumber\\[-8pt]\\[-8pt]
&&\qquad\hspace*{0pt}\leq e^\ep1_{\{\|h_n\|\leq M\}}\int_{D({\theta}_n,\rho_n)}
\prod_{i=1}^n\frac{p_{{\theta}_n,\eta}}{p_{{\theta}_0,\eta
_0}}(X_i) \,d\Pi
_H(\eta)
\Biggr) > 1-\delta.
\nonumber
\end{eqnarray}
We continue with the integral over $D({\theta}_n,\rho_n)$ under
the restriction $\|h_n\|\leq M$ and parametrize the model
locally in terms of $({\theta},\zeta)$ [see~(\ref{eqrepara})]
%
%e4.6 #&#
%
\begin{equation}
\label{eqrepint}
\int_{D({\theta}_n,\rho_n)}
\prod_{i=1}^n\frac{p_{{\theta}_n,\eta}}{p_{{\theta}_0,\eta_0}}(X_i)
\,d\Pi_H(\eta)\\
= \int_{B(\rho_n)}
\prod_{i=1}^n\frac{q_{{\theta}_n,\zeta}}{q_{{\theta}_0,0}}(X_i)
\,d\Pi( \zeta\bigmv {\theta}={\theta}_n),\hspace*{-30pt}
\end{equation}
where $\Pi( \cdot| {\theta} )$ denotes the prior for $\zeta$ given
${\theta}$, that is, $\Pi_H$ translated over~$\eta^*({\theta})$.
Next we note
that by Fubini's theorem and the domination condition~(\ref{eqnormdom}),
there exists a constant $L>0$ such that
\begin{eqnarray*}
&&\Biggl| P_0^n\int_{B(\rho_n)}
\prod_{i=1}^n\frac{q_{{\theta}_n,\zeta}}{q_{{\theta}_0,0}}(X_i)
\bigl(d\Pi(\zeta\bigmv {\theta}={\theta}_n)
-d\Pi(\zeta\bigmv {\theta}={\theta}_0)\bigr)
\Biggr|\\
&&\qquad\leq L
\bigl| \Pi\bigl( B(\rho_n)\bigmv {\theta}={\theta}_n \bigr)
-\Pi\bigl( B(\rho_n)\bigmv {\theta}={\theta}_0 \bigr) \bigr|
\end{eqnarray*}
for large enough $n$. Since the least-favorable submodel is
stochastically LAN, Lemma~\ref{lemtranslate} asserts
that the difference on the r.h.s. of the above display is~$o(1)$,
so that
%
%e4.7 #&#
%
\begin{eqnarray}
\label{eqshiftprior}
&&
\int_{B(\rho_n)}
\prod_{i=1}^n\frac{q_{{\theta}_n,\zeta}}{q_{{\theta}_0,0}}(X_i)
\,d\Pi(\zeta\bigmv {\theta}={\theta}_n)\nonumber\\[-8pt]\\[-8pt]
&&\qquad= \int_{B(\rho_n)}
\prod_{i=1}^n\frac{q_{{\theta}_n,\zeta}}{q_{{\theta}_0,0}}(X_i)
\,d\Pi(\zeta) + o_{P_0}(1),\nonumber
\end{eqnarray}
where we use the notation $\Pi(A)=\Pi( \zeta\in A | {\theta
}={\theta}_0 )$
for brevity. We define for all $\zeta$, $\ep>0$, $n\geq1$ the events
$F_n(\zeta,\ep)=\{\sup_h|G_n(h,\zeta)-G_n(h,0)|
\leq\ep\}$. With~(\ref{eqnormdom}) as a domination condition,
Fatou's lemma and the fact that $F_n^c(0,\ep)=\varnothing$ lead to
%
%e4.8 #&#
%
\begin{eqnarray}
\label{eqintQF}
&&
\limsup_{n\rightarrow\infty}
\int_{B(\rho_n)} Q^n_{{\theta}_n,\zeta}(F_n^c(\zeta,\ep
))
\,d\Pi(\zeta)\nonumber\\[-8pt]\\[-8pt]
&&\qquad\leq\int\limsup_{n\rightarrow\infty} 1_{B(\rho_n)\setminus\{0\}
}(\zeta)
Q^n_{{\theta}_n,\zeta}(F_n^c(\zeta,\ep)) \,d\Pi(\zeta)
=0\nonumber
\end{eqnarray}
[again using~(\ref{eqnormdom}) in the last step]. Combined with
Fubini's theorem, this suffices to conclude that
%
%e4.9 #&#
%
\begin{equation}
\label{eqintrobn}
\int_{B(\rho_n)}
\prod_{i=1}^n\frac{q_{{\theta}_n,\zeta}}{q_{{\theta}_0,0}}(X_i)
\,d\Pi(\zeta)
= \int_{B(\rho_n)}
\prod_{i=1}^n\frac{q_{{\theta}_n,\zeta}}{q_{{\theta}_0,0}}(X_i)
1_{F_n(\zeta,\ep)} \,d\Pi(\zeta) + o_{P_0}(1),\hspace*{-32pt}
\end{equation}
and we continue with the first term on the right-hand side. By stochastic
local asymptotic normality for every $\zeta$,
expansion~(\ref{eqqlan}) of the log-likelihood implies that
%
%e4.10 #&#
%
\begin{equation}
\label{eqLANintegrand}
\prod_{i=1}^n \frac{q_{{\theta}_n,\zeta}}{q_{{\theta}_0,0}}(X_i)
=\prod_{i=1}^n \frac{q_{{\theta}_0,\zeta}}{q_{{\theta}_0,0}}(X_i)
e^{G_n(h_n,\zeta)+R_n(h_n,\zeta)},
\end{equation}
where the rest term is of order $o_{Q_{{\theta}_0,\zeta}}(1)$. Accordingly,
we define, for every $\zeta$, the events $A_n(\zeta,\ep)=\{
|R_n(h_n,\zeta)|\leq\ft12\ep\}$, so that
$Q^n_{{\theta}_0,\zeta}(A_n^c(\zeta,\ep))\rightarrow0$. Contiguity
then implies that $Q^n_{{\theta}_n,\zeta}(A^c_n(\zeta,\ep
))\rightarrow0$
as well. Reasoning as in~(\ref{eqintrobn}) we see that
%
%e4.11 #&#
%
\begin{eqnarray}
\label{eqintroan}
&&
\int_{B(\rho_n)} \prod_{i=1}^n\frac{q_{{\theta}_n,\zeta}}
{q_{{\theta}_0,0}}(X_i) 1_{F_n(\zeta,\ep)}
\,d\Pi(\zeta)\nonumber\\[-8pt]\\[-8pt]
&&\qquad=
\int_{B(\rho_n)}\prod_{i=1}^n\frac{q_{{\theta}_n,\zeta}}
{q_{{\theta}_0,0}}(X_i) 1_{A_n(\zeta,\ep)\cap F_n(\zeta,\ep)}
\,d\Pi(\zeta)
+ o_{P_0}(1).
\nonumber
\end{eqnarray}
For fixed $n$ and $\zeta$ and for all
$(X_1,\ldots,X_n)\in A_n(\zeta,\ep)\cap F_n(\zeta,\ep)$,
\[
\Biggl| \log\prod_{i=1}^n\frac{q_{{\theta}_n,\zeta}}
{q_{{\theta}_0,0}}(X_i) - G_n(h_n,0) \Biggr|
\leq2\ep,
\]
so that the first term on the right-hand side of~(\ref{eqintroan}) satisfies
the bounds
%
%e4.12 #&#
%
\begin{eqnarray}
\label{eqexpcorrection}
&&e^{G_n(h_n,0)-2\ep}\int_{B(\rho_n)}
\prod_{i=1}^n\frac{q_{{\theta}_0,\zeta}}{q_{{\theta}_0,0}}(X_i)
1_{A_n(\zeta,\ep)\cap F_n(\zeta,\ep)}
\,d\Pi(\zeta)\nonumber\\
&&\qquad \leq
\int_{B(\rho_n)}\prod_{i=1}^n\frac{q_{{\theta}_n,\zeta}}
{q_{{\theta}_0,0}}(X_i)
1_{A_n(\zeta,\ep)\cap F_n(\zeta,\ep)}
\,d\Pi(\zeta)\\
&&\qquad\leq
e^{G_n(h_n,0)+2\ep}\int_{B(\rho_n)}
\prod_{i=1}^n\frac{q_{{\theta}_0,\zeta}}{q_{{\theta}_0,0}}(X_i)
1_{A_n(\zeta,\ep)\cap F_n(\zeta,\ep)}
\,d\Pi(\zeta).
\nonumber
\end{eqnarray}
The integral factored into lower and upper bounds can be
relieved of the indicator for $A_n\cap F_n$ by reversing the
argument that led to~(\ref{eqintrobn}) and~(\ref{eqintroan})
(with~${\theta}_0$ replacing ${\theta}_n$), at the expense of an
$e^{o_{P_0}(1)}$-factor. Substituting in~(\ref{eqexpcorrection})
and using, consecutively,~(\ref{eqintroan}),~(\ref{eqintrobn}),
(\ref{eqshiftprior}) and~(\ref{eqintD}) for the bounded
integral, we find
\[
e^{G_n(h_n,0)-3\ep+o_{P_0}(1)} s_n(0)
\leq s_n(h_n) \leq
e^{G_n(h_n,0)+3\ep+o_{P_0}(1)}s_n(0).
\]
Since this holds with arbitrarily small $0<\ep'<\ep$ for large enough
$n$, it proves~(\ref{eqilan}).
\end{pf}

With regard to the nuisance rate $(\rho_n)$, we first note that our
proof of Theorem~\ref{thmsbvmone} fails if the slowest rate required
to satisfy~(\ref{eqnormdom}) vanishes \textit{faster} then the optimal
rate for convergence under $n^{-1/2}$-perturbation [as determined in
(\ref{eqminimaxrate}) and~(\ref{eqsuffprior})].

However, the rate $(\rho_n)$ does not appear in assertion~(\ref{eqilan}),
so if said contradiction between conditions~(\ref{eqnormdom}) and
(\ref{eqminimaxrate})/(\ref{eqsuffprior}) do not occur, the
sequence~$(\rho_n)$ can remain entirely internal to the proof of
Theorem~\ref{thmilanone}. More particularly, if condition~(\ref{eqnormdom})
holds for \textit{any} $(\rho_n)$ such that $n\rho_n^2\rightarrow\infty$,
integral LAN only requires consistency under $n^{-1/2}$-perturbation
at \textit{some} such~$(\rho_n)$. In that case, we may appeal to
Corollary~\ref{corconspert} instead of Theorem~\ref{thmpertroc},
thus relaxing conditions on model entropy and nuisance prior.
The following lemma shows that a first-order Taylor expansion of
likelihood ratios combined with a boundedness condition on certain
Fisher information coefficients is enough to enable use of
Corollary~\ref{corconspert} instead of Theorem~\ref{thmpertroc}.

%le4.3 #&#
%
\begin{lemma}
\label{lemUdom}
Let $\Tht$ be one-dimensional. Assume that there exists
a $\rho>0$ such that for every $\zeta\in B(\rho)$ and all $x$ in the
samplespace, the map ${\theta}\mapsto\log(q_{{\theta},\zeta
}/q_{{\theta}_0,\zeta
})(x)$ is
continuously differentiable on $[{\theta}_0-\rho,{\theta}_0+\rho]$ with
Lebes\-gue-integrable derivative $g_{{\theta},\zeta}(x)$ such that
%
%e4.13 #&#
%
\begin{equation}
\label{eqUdomcond}
\sup_{\zeta\in B(\rho)} \sup_{\{{\theta}\dvtx|{\theta}-{\theta
}_0|<\rho\}}
Q_{{\theta},\zeta}g_{{\theta},\zeta}^2<\infty.
\end{equation}
Then, for every $\rho_n\downarrow0$
and all bounded, stochastic $(h_n)$, $U_n(\rho_n,h_n)=O(1)$.
\end{lemma}
\begin{pf}
Let $(h_n)$ be stochastic and upper-bounded by $M>0$. For every~$\zeta$ and all $n\geq1$,
\begin{eqnarray*}
Q_{{\theta}_0,\zeta}^n\Biggl|\prod_{i=1}^n\frac{q_{{\theta
}_n(h_n),\zeta}}
{q_{{\theta}_0,\zeta}}(X_i)-1\Biggr|
&=& Q_{{\theta}_0,\zeta}^n\Biggl|
\int_{{\theta}_0}^{{\theta}_n(h_n)}\sum_{i=1}^ng_{{\theta}',\zeta}(X_i)
\prod_{j=1}^n\frac{q_{{\theta}',\zeta}}{q_{{\theta}_0,\zeta
}}(X_j) \,d{\theta}
'\Biggr|\\
&\leq&\int_{{\theta}_0-{M}/{\sqrt{n}}}^{{\theta}_0+
{M}/{\sqrt{n}}}
Q_{{\theta}',\zeta}^n\Biggl|\sum_{i=1}^ng_{{\theta}',\zeta
}(X_i)\Biggr| \,d{\theta}'\\
&\leq&\sqrt{n}\int_{{\theta}_0-{M}/{\sqrt{n}}}^{{\theta
}_0+{M}/{\sqrt{n}}}
\sqrt{Q_{{\theta}',\zeta}g_{{\theta}',\zeta}^2} \,d{\theta}',
\end{eqnarray*}
where the last step follows from the Cauchy--Schwarz inequality. For
large enough $n$, $\rho_n<\rho$ and the square-root of~(\ref{eqUdomcond})
dominates the difference between $U(\rho,h_n)$ and $1$.
\end{pf}

%%%%%%%%%%%%%%%%%%%%%%%%%%%%%%%%%%%%%%%%%%%%%%%%%%%%%%%%%%%%%%%%%%%%%%%%%%%%%%%%%%

%s5 #&#
\section{Posterior asymptotic normality}
\label{secpan}

Under the assumptions formulated before
Theorem~\ref{thmsbvmone}, the marginal posterior density
$\pi_n(\cdot|X_1,\ldots,X_n)\dvtx\Tht\rightarrow\RR$ for the
parameter of interest with respect to the prior $\Pi_\Tht$
equals
%
%e5.1 #&#
%
\begin{equation}
\label{eqDefPD}
\pi_n({\theta}|X_1,\ldots,X_n) = {S_n({\theta})}\Big/
{ {\int_\Tht S_n({\theta}') \,d\Pi_{\Tht}({\theta}')}},
\end{equation}
$P_0^n$-almost-surely. One notes that this form is equal to that
of a \textit{parametric} posterior density, but with the
parametric likelihood replaced by the integrated likelihood
$S_n$. By implication, the proof of the parametric Bernstein--von Mises
theorem can be applied to its semiparametric generalization, if we
impose sufficient conditions for the parametric likelihood on $S_n$
instead. Concretely, we replace the smoothness requirement for the
likelihood in Theorem~\ref{thmparaBvM} by~(\ref{eqilan}). Together
with a condition expressing marginal posterior convergence at
parametric rate,~(\ref{eqilan}) is sufficient to derive
asymptotic normality of the posterior; cf.~(\ref{eqassertBvM}).
%This shortcut is illustrated further by the following perspective.
%For given $\tht$ and $n$, $s_n(n^{1/2}(\tht-\tht_0))$ is a
%probability density for the stochastic vector $(X_1,\ldots,X_n)$
%with respect to $P_0^n$, corresponding to the $\tht$-conditioned
%($\Pi_H$-prior predictive) distribution,
% \tilde{P}_{n,\tht}(B) = P_0^n(1_B
% s_n(\sqrt{n}(\tht-\tht_0))),
%(where $B$ measurable in the $n$-fold product of the
%samplespace). Indeed, keeping $n$ fixed, we may view the map
%$\tht\mapsto\tilde{P}_{n,\tht}$ as a parametric model with a
%prior $\Pi_\Tht$ that is thick at $\tht_0$. Condition
%(\ref{eqilan}) then amounts to
%stochastic local asymptotic normality of this parametric model
%and condition {\it{(iv)}} of theorem~\ref{thmsbvmone} to parametric
%rate-optimality of its posterior. This conceptual simplification
%comes at a price, though: firstly, this parametric model is
%misspecified, \ie\ there is no $\tht\in\Tht$ such that
%$P^n_0=\tilde{P}_{n,\tht}$. Secondly, although we have assumed
%that the sample is distributed \iid, in the parametric model
%above $X_1,\ldots, X_n$ are \textit{not} independent, instead the
%sample $(X_1,\ldots, X_n)$ satisfies the weaker property of
%exchangeability under $\tilde{P}_{n,\tht}$ for every $\tht$,
%in accordance with De Finetti's theorem. Although this enables
%application of methods put forth in Kleijn and van der Vaart
%%\cite{Kleijn07}, in the present context, results are sharper
%if we take into account the semiparametric background of the
%quantities $s_n(h)$.
%
%th5.1 #&#
%
\begin{Thm}[(Posterior asymptotic normality)]
\label{thmpan}
Let $\Tht$ be open in $\RR^k$ with a prior $\Pi_\Tht$ that is
thick at ${\theta}_0$. Suppose that for large enough $n$, the map
$h\mapsto s_n(h)$ is continuous $P_0^n$-almost-surely. Assume
that there exists an $L_2(P_0)$-function $\effscore_{{\theta}_0,\eta_0}$
such that for every $(h_n)$ that is bounded in probability,~(\ref{eqilan})
holds, $P_0\effscore_{{\theta}_0,\eta_0}=0$ and
$\effFI_{{\theta}_0,\eta_0}$ is nonsingular. Furthermore suppose that
for every $(M_n)$, $M_n\rightarrow\infty$,
we have
%
%e5.2 #&#
%
\begin{equation}
\label{eqsqrtn}
\Pi_n( \|h\|\leq M_n\bigmv X_1,\ldots,X_n )\convprob
{P_0} 1.
\end{equation}
Then the sequence of marginal posteriors for ${\theta}$ converges
to a normal distribution in total variation,
\[
\sup_{A}\bigl|
\Pi_n( h\in A\bigmv X_1,\ldots,X_n )
- N_{\effDelta_n,\effFI_{{\theta}_0,\eta_0}^{-1}}(A)
\bigr| \convprob{P_0} 0,
\]
centered on $\effDelta_n$ with covariance matrix $\effFI_{{\theta
}_0,\eta_0}^{-1}$.
\end{Thm}
\begin{pf}
The proof is identical to that of Theorem 2.1 in~\cite{Kleijn07}
upon replacement of parametric likelihoods with integrated likelihoods.
\end{pf}

There is room for relaxation of the requirements
on model entropy and minimal prior mass, if the limit~(\ref{eqnormdom})
holds in a fixed neighborhood of $\eta_0$. The following
corollary applies whenever~(\ref{eqnormdom}) holds for \textit{any rate}
$(\rho_n)$. The simplifications are such that the entropy and
prior mass conditions become comparable to those for
Schwartz's posterior consistency theorem~\cite{Schwartz65}, rather
than those for posterior rates of convergence following Ghosal,
Ghosh and van der Vaart~\cite{Ghosal00}.
%
%co5.2 #&#
%
\begin{corollary}[(Semiparametric Bernstein--von Mises,
rate-free)] \label{corsimplesbvm}
Let $X_1$, $X_2,\ldots$ be i.i.d.-$P_0$,
with $P_0\in\scrP$, and let $\Pi_\Tht$ be thick at ${\theta}_0$.
Suppose that for large enough~$n$, the map $h\mapsto s_n(h)$ is
continuous $P_0^n$-almost-surely. Also assume that ${\theta}\mapsto
Q_{{\theta},\zeta}$ is stochastically LAN in the ${\theta}$-direction,
for all $\zeta$ in an $r_H$-neighborhood of $\zeta=0$ and that the
efficient Fisher information $\effFI_{{\theta }_0.\eta_0}$ is
nonsingular. Furthermore, assume that:
\begin{longlist}[(iii)]
\item[(i)] For all $\rho>0$, the Hellinger metric entropy satisfies,
$N(\rho,H,d_H) < \infty$ and the nuisance prior satisfies
$\Pi_H( K(\rho) ) > 0$.
\item[(ii)] For every $M>0$, there exists an $L>0$ such that for all
$\rho>0$ and large enough $n$, $K(\rho) \subset K_n(L\rho,M)$.

Assume also that for every bounded, stochastic $(h_n)$:
\item[(iii)] There exists an $r>0$ such that, $U_n(r,h_n)=O(1)$.
\item[(iv)] Hellinger distances satisfy,
$\sup_{\eta\in H}H(P_{{\theta}_n(h_n),\eta},P_{{\theta}_0,\eta
})=O(n^{-1/2})$,

and that
\item[(v)] For every $(M_n)$, $M_n\rightarrow\infty$, the posterior
satisfies,
\[
\Pi_n( \|h\|\leq M_n\bigmv X_1,\ldots,X_n )\convprob
{P_0} 1.
\]
\end{longlist}
Then the sequence of marginal posteriors for ${\theta}$ converges in
total variation to a normal distribution,
\[
\sup_{A}\bigl|
\Pi_n( h\in A\bigmv X_1,\ldots,X_n )
- N_{\effDelta_n,\effFI_{{\theta}_0,\eta_0}^{-1}}(A)
\bigr| \convprob{P_0} 0,
\]
centered on $\effDelta_n$ with covariance matrix $\effFI_{{\theta
}_0,\eta_0}^{-1}$.
\end{corollary}
\begin{pf}
Under conditions (i), (ii), (iv) and the stochastic
LAN assumption, the assertion of Corollary~\ref{corconspert} holds. Due
to condition (iii), condition~(\ref{eqnormdom}) is satisfied
for large enough $n$. Condition (v) then suffices for the assertion
of Theorem~\ref{thmpan}.
\end{pf}

A critical note can be made regarding the qualification ``rate-free'' of
Corollary~\ref{corsimplesbvm}: although the nuisance rate does not make
an explicit appearance, rate restrictions may arise upon further analysis
of condition~(v). Indeed this is the case in the example of
Section~\ref{secplr}, where smoothness requirements on the regression
family are interpretable as restrictions on the nuisance rate.
However, semiparametric models exist, in which no restrictions on
nuisance rates arise in this way: if $H$ is a convex
subspace of a~linear space, and the dependence $\eta\mapsto P_{{\theta
},\eta}$
is linear (a so-called \textit{convex-linear} model, e.g., mixture models,
errors-in-variables regression and other information-loss models), the
construction of suitable tests (cf. Le Cam~\cite{LeCam86},
Birg\'e~\cite{Birge83,Birge84}) does not involve Hellinger
metric entropy numbers or restrictions on nuisance rates of convergence.
Consequently there exists a class of semiparametric examples for which
Corollary~\ref{corsimplesbvm} stays rate-free even after further
analysis of its condition (v).

As shown in~\cite{Kleijn07}, the particular form of the limiting
posterior in Theorem~\ref{thmpan} is a consequence of local asymptotic
normality, in this case imposed through~(\ref{eqilan}). The
marginal posterior converges exactly to the asymptotic sampling
distribution of a frequentist best-regular estimator as a consequence.
Other expansions (e.g., in LAN models for non-i.i.d. data
or under the condition of \textit{local asymptotic exponentiality}
(Ibragimov and Has'mins\-kii~\cite{Ibragimov81})) can be dealt
with in the same manner if we adapt the limiting form of the posterior
accordingly, giving rise to other (e.g., one-sided exponential)
limit distributions (see Kleijn and Knapik~\cite{Kleijn10a}).

%%%%%%%%%%%%%%%%%%%%%%%%%%%%%%%%%%%%%%%%%%%%%%%%%%%%%%%%%%%%%%%%%%%%%%%%%%%%%%%%%%

%s6 #&#
\section{Marginal posterior convergence at parametric rate}
\label{secmarg}

Condition~(\ref{eqsqrtn}) in Theorem~\ref{thmpan} requires that the
posterior measures of a sequence of model subsets of the form
%
%e6.1 #&#
%
\begin{equation}
\label{eqstrip}
\Tht_n\times H = \bigl\{ ({\theta},\eta)\in\Tht\times H
\dvtx\sqrt{n}\|{\theta}-{\theta}_0\|\leq M_n\bigr\}
\end{equation}
converge to one in $P_0$-probability, for every sequence
$(M_n)$ such that $M_n\rightarrow\infty$. Essentially, this
condition enables us to restrict the proof of Theorem~\ref{thmpan}
to the shrinking domain in which~(\ref{eqilan}) applies.
%Marginal posteriors have not received much specific attention
%in the literature on posterior asymptotics thus far.
%Questions concerning testing in the presence of nuisance
%parameters~\cite{Choi96,Bickel04} lie at the centre of
%this problem.
In this section, we consider two distinct
approaches: the first (Lemma~\ref{lemlehmann}) is based on
bounded likelihood ratios (see also condition (B3) of Theorem 8.2 in
Lehmann and Casella~\cite{Lehmann98}). The second is
based on the behavior of misspecified parametric posteriors
(Theorem~\ref{thmrocBayes}). The latter construction illustrates
the intricacy of this section's subject most clearly and provides
some general insight. Methods proposed here are neither compelling
nor exhaustive; we simply put forth several possible approaches
and demonstrate the usefulness of one of them in Section~\ref{secplr}.
%
%le6.1 #&#
%
\begin{lemma}[{[Marginal parametric rate (I)]}]
\label{lemlehmann}
Let the sequence of maps ${\theta}\mapsto S_n({\theta})$ be
$P_0$-almost-surely
continuous and such that~(\ref{eqilan}) is satisfied. Furthermore,
assume that there exists a constant $C>0$ such that for any~$(M_n)$, $M_n\rightarrow\infty$,
%
%e6.2 #&#
%
\begin{equation}
\label{eqlehmann}
P_0^n\biggl( \sup_{\eta\in H} \sup_{{\theta}\in\Tht^c_n}
\PP_n\log\frac{p_{{\theta},\eta}}{p_{{\theta}_0,\eta}}
\leq-\frac{C M_n^2}{n} \biggr)
\rightarrow1.
\end{equation}
Then, for any nuisance prior $\Pi_H$ and parametric prior
$\Pi_{\Tht}$, thick at ${\theta}_0$,
%
%e6.3 #&#
%
\begin{equation}
\label{eqrocULR}
\Pi( n^{1/2}\|{\theta}-{\theta}_0\|>M_n \bigmv  X_1,\ldots
,X_n )
\convprob{P_0} 0
\end{equation}
for any $(M_n)$, $M_n\rightarrow\infty$.
\end{lemma}
\begin{pf}
Let $(M_n)$, $M_n\rightarrow\infty$ be given. Define $(A_n)$ to be the
events in~(\ref{eqlehmann}) so that $P_0^n(A_n^c)=o(1)$ by assumption.
In addition, let
\[
B_n = \biggl\{\int_\Tht S_n({\theta}) \,d\Pi_\Tht({\theta})
\geq e^{-C M_n^2/2} S_n({\theta}_0) \biggr\}.
\]
By~(\ref{eqilan}) and Lemma~\ref{lemdenom}, $P_0^n(B_n^c)=o(1)$ as well.
Then
\begin{eqnarray*}
&&P_0^n\Pi({\theta}\in\Tht^c_n\mid X_1,\ldots,X_n)\\
&&\quad\leq P_0^n\Pi({\theta}\in\Tht^c_n\mid X_1,\ldots,X_n)
1_{A_n\cap B_n} + o(1)\\
&&\quad\leq e^{C M_n^2/2} P_0^n\Biggl( S_n({\theta}_0)^{-1}
\int_H\int_{\Tht^c_n}
\prod_{i=1}^n\frac{p_{{\theta},\eta}}{p_{{\theta}_0,\eta}}(X_i)
\prod_{i=1}^n\frac{p_{{\theta}_0,\eta}}{p_{{\theta}_0,\eta
_0}}(X_i)
\,d\Pi_\Tht \,d\Pi_H 1_{A_n} \Biggr)\\
&&\quad\quad{} +o(1)\\
&&\quad = o(1),
\end{eqnarray*}
which proves~(\ref{eqrocULR}).
\end{pf}

Although applicable directly in the model of Section~\ref{secplr},
most other examples would require variations.
Particularly, if the full, nonparametric posterior is known to
concentrate on a sequence
of model subsets $(V_n)$, then Lemma~\ref{lemlehmann} can be preceded
by a decomposition of $\Tht\times H$ over $V_n$ and $V_n^c$, reducing
condition~(\ref{eqlehmann}) to a supremum over $V_n^c$ (see Section 2.4
in Kleijn~\cite{Kleijn03} and the discussion following the
following theorem).

Our second approach assumes such concentration of the
posterior on model subsets, for example, deriving from nonparametric
consistency in a~suitable form. Though the proof of
Theorem~\ref{thmrocBayes} is rather straightforward, combination
with results in misspecified parametric models~\cite{Kleijn07} leads
to the observation that marginal parametric rates of convergence can
be ruined by a bias.
%
%th6.2 #&#
%
\begin{Thm}[{[Marginal parametric rate (II)]}]
\label{thmrocBayes}
Let $\Pi_\Tht$ and $\Pi_H$ be given. Assume that there exists
a sequence $(H_n)$ of subsets of $H$, such that the following
two conditions hold:
\begin{longlist}[(ii)]
\item[(i)] The nuisance posterior concentrates on $H_n$ asymptotically,
%
%e6.4 #&#
%
\begin{equation}
\label{eqetacons}
\Pi( \eta\in H\setminus H_n\bigmv X_1,\ldots,X_n
)\convprob{P_0}0.
\end{equation}
\item[(ii)] For every $(M_n)$, $M_n\rightarrow\infty$,
%
%e6.5 #&#
%
\begin{equation}
\label{eqptwise}
P_0^n\sup_{\eta\in H_n}\Pi( n^{1/2}\|{\theta}-{\theta}_0\|> M_n
\bigmv \eta,X_1,\ldots,X_n ) \rightarrow0.
\end{equation}
\end{longlist}
Then the marginal posterior for ${\theta}$ concentrates
at parametric rate, that is,
\[
\Pi( n^{1/2}\|{\theta}-{\theta}_0\|> M_n\bigmv \eta
,X_1,\ldots,X_n)
\convprob{P_0}0
\]
for every sequence $(M_n)$, $M_n\rightarrow\infty$.\vadjust{\goodbreak}
\end{Thm}
\begin{pf} Let $(M_n)$, $M_n\rightarrow\infty$ be given, and
consider the posterior for the complement of~(\ref{eqstrip}). By
assumption (i) of the theorem and Fubini's theorem,
\begin{eqnarray*}
&&P_0^n\Pi( {\theta}\in\Tht_n^c\bigmv X_1,\ldots,X_n)\\
&&\qquad\leq P_0^n\int_{H_n}
\Pi( {\theta}\in\Tht_n^c\bigmv \eta,X_1,\ldots,X_n)
\,d\Pi( \eta\bigmv X_1,\ldots,X_n)
+ o(1)\\
&&\qquad\leq P_0^n\sup_{\eta\in H_n}\Pi( n^{1/2}\|{\theta}-{\theta
}_0\|> M_n
\bigmv \eta,X_1,\ldots,X_n) + o(1),
\end{eqnarray*}
the first term of which is $o(1)$ by assumption (ii) of the theorem.
\end{pf}

Condition (ii) of Theorem~\ref{thmrocBayes} has an
interpretation in terms of misspecified parametric models
(Kleijn and van der Vaart~\cite{Kleijn07} and Kleijn
\cite{Kleijn03}). For fixed $\eta\in H$, the $\eta$-conditioned
posterior on the parametric model $\scrP_\eta=\{P_{{\theta},\eta
}\dvtx{\theta}\in
\Tht\}$
is required to concentrate in $n^{-1/2}$-neighborhoods
of ${\theta}_0$ under~$P_0$. However, this misspecified posterior
concentrates around $\Tht^*(\eta)\subset\Tht$, the set of points
in $\Tht$ where the Kullback--Leibler divergence of $P_{{\theta},\eta}$
with respect to $P_0$, is minimal. Assuming that $\Tht^*(\eta)$
consists of a unique minimizer~${\theta}^*(\eta)$, the dependence
of the Kullback--Leibler divergence on $\eta$ must be such that
%
%e6.6 #&#
%
\begin{equation}
\label{eqnearstraight}
\sup_{\eta\in H_n} \|{\theta}^*(\eta)-{\theta}_0\| = o
(n^{-1/2})
\end{equation}
in order for posterior concentration to occur on the strips
(\ref{eqstrip}). In other words, minimal Kullback--Leibler
divergence may bias the (points of convergence of) $\eta$-conditioned
parametric posteriors to such an extent that consistency of the
marginal posterior for ${\theta}$ is ruined.

The occurrence of this bias is a property of the semiparametric
model rather than a peculiarity of the Bayesian approach: when
(point-)estimating with solutions to score equations, for example,
the same bias occurs (see, e.g., Theorem~25.59 in~\cite{vdVaart98}
and subsequent discussion). Frequentist literature also offers
some guidance toward mitigation of this circumstance. First of
all, it is noted that the bias indicates the existence of a better
(i.e., bias-less) choice of parametrization to ask the relevant
semiparametric question. If the parametrization is fixed,
alternative point-estimation methods may resolve bias, for
example, through replacement of score equations by general estimating
equations (see, e.g., Section 25.9 in~\cite{vdVaart98}),
loosely equivalent to introducing a suitable penalty in a
likelihood maximization procedure.

For a so-called \textit{curve-alignment model} with Gaussian prior,
the no-bias problem has been addressed and resolved in a fully
Bayesian manner by Castillo~\cite{Castillo11}: like a
penalty in an ML procedure, Castillo's (rather subtle choice of)
prior guides the procedure away from the biased directions and
produces Bernstein--von Mises efficiency of the marginal posterior.
A most interesting question concerns generalization of
Castillo's intricate construction to more general Bayesian context.

Recalling definitions~(\ref{eqdefSn}) and~(\ref{eqdefsn}), we
conclude this section with a~lemma used in the proof of
Lemma~\ref{lemlehmann} to lower-bound the denominator
of the marginal posterior.
%
%le6.3 #&#
%
\begin{lemma}
\label{lemdenom}
Let the sequence of maps ${\theta}\mapsto S_n({\theta})$ be
$P_0$-almost-surely
continuous and such that~(\ref{eqilan}) is satisfied. Assume that
$\Pi_{\Tht}$ is thick at ${\theta}_0$ and denoted by $\Pi_n$ in the local
parametrization in terms of $h$. Then
%
%e6.7 #&#
%
\begin{equation}
\label{eqdenom}
P_0^n\biggl( \int s_n(h) \,d\Pi_n(h)<a_n s_n(0) \biggr)\rightarrow0
\end{equation}
for every sequence $(a_n)$, $a_n\downarrow0$.
\end{lemma}
\begin{pf} Let $M>0$ be given, and define $C=\{h\dvtx\|h\|\leq M\}$.
Denote the rest-term in~(\ref{eqilan}) by $h\mapsto R_n(h)$.
By continuity of ${\theta}\mapsto S_n({\theta})$, ${\sup_{h\in C}}|R_n(h)|$
converges to zero in $P_0$-probability. If we choose a sequence
$(\kappa_n)$ that converges to zero slowly enough, the corresponding
events $B_n = \{ {\sup_C} |R_n(h)|
\leq\kappa_n \}$, satisfy $P_0^n(B_n)\rightarrow1$. Next,
let $(K_n)$, $K_n\rightarrow\infty$ be given. There exists a~$\pi>0$
such that $\inf_{h\in C} d\Pi_n/d\mu(h) \geq\pi$, for large enough
$n$. Combining, we find
%
%e6.8 #&#
%
\begin{eqnarray}
\label{eqdenomone}
&&
P_0^n\biggl( \int\frac{s_n(h)}{s_n(0)} \,d\Pi_n(h)\leq e^{-K_n^2}
\biggr)\nonumber\\[-8pt]\\[-8pt]
&&\qquad\leq
P_0^n\biggl( \biggl\{ \int_C
\frac{s_n(h)}{s_n(0)} \,d\mu(h)\leq\pi^{-1} e^{-K_n^2}\biggr\}
\cap B_n \biggr) + o(1).
\nonumber
\end{eqnarray}
On $B_n$, the integral LAN expansion is lower
bounded so that, for large enough $n$,
%
%e6.9 #&#
%
\begin{eqnarray}
\label{eqdenomtwo}
&&P_0^n\biggl( \biggl\{ \int_C
\frac{s_n(h)}{s_n(0)}\, d\mu(h)\leq\pi^{-1} e^{-K_n^2}\biggr\}
\cap B_n \biggr)\nonumber\\[-8pt]\\[-8pt]
&&\qquad\leq
P_0^n\biggl( \int_C e^{h^T\GG_n\effscore_{{\theta}_0,\eta_0}}
\,d\mu(h)\leq\pi^{-1}e^{-K_n^2/4}
\biggr)\nonumber
\end{eqnarray}
since $\kappa_n\leq\ft12K_n^2$ and ${\sup_{h\in C}}|h^T\effFI
_{{\theta}_0,\eta_0}h|
\leq M^2 \|\effFI_{{\theta}_0,\eta_0}\|\leq\ft14K_n^2$, for large\break
enough~$n$.
Conditioning $\mu$ on $C$, we apply Jensen's inequality to note that,
for large enough~$n$,
\begin{eqnarray*}
&&P_0^n\biggl( \int_C e^{h^T\GG_n\effscore_{{\theta}_0,\eta_0}}
\, d\mu(h)\leq\pi^{-1}e^{-K_n^2/4} \biggr)\\
&&\qquad\leq
P_0^n\biggl( \int h^T\GG_n\effscore_{{\theta}_0,\eta_0}
\,d\mu(h|C) \leq-\ft18K_n^2 \biggr)
\end{eqnarray*}
since\vspace*{1pt} $-{\log\pi}\mu(C)\leq\ft18K_n^2$, for large enough $n$.
The probability on the right is bounded further by Chebyshev's
and Jensen's inequalities and can be shown to be of order $O(K_n^{-4})$.
Combining with~(\ref{eqdenomone}) and~(\ref{eqdenomtwo}) then
proves~(\ref{eqdenom}).
\end{pf}

%%%%%%%%%%%%%%%%%%%%%%%%%%%%%%%%%%%%%%%%%%%%%%%%%%%%%%%%%%%%%%%%%%%%%%%%%%%%%%%%%%

%s7 #&#
\section{Semiparametric regression}
\label{secplr}

The \textit{partial linear regression} model describes the observation of an
i.i.d. sample $X_1,X_2,\ldots$ of triplets $X_i=(U_i,V_i,\allowbreak Y_i)\in\RR^3$,
each assumed to be related through the regression equation
%
%e7.1 #&#
%
\begin{equation}
\label{eqplrmodel}
Y = {\theta}_0 U + \eta_0(V) + e,
\end{equation}
where $e\sim N(0,1)$ is independent of $(U,V)$. Interpreting $\eta_0$
as a nuisance parameter, we wish to estimate ${\theta}_0$. It is assumed
that $(U,V)$ has an unknown distribution~$P$, Lebesgue absolutely
continuous with density $p\dvtx\RR^2\rightarrow\RR$. The distribution $P$
is assumed to be such that $PU=0$, $PU^2=1$ and $PU^4<\infty$. At a
later stage, we also impose $P(U-\mathrm{E}[U|V])^2>0$ and a smoothness
condition on the conditional expectation $v\mapsto\mathrm{E}[U|V=v]$.

As is well known~\cite{Chen91,Bickel98,Mammen97,vdVaart98}, penalized
ML estimation in a smoothness class of regression functions leads to
a consistent estimate of the nuisance and efficient point-estimation
of the parameter of interest. The necessity of a penalty signals that
the choice of a prior for the nuisance is a critical one.
Kimeldorf and Wahba~\cite{Kimeldorf70} assume that the regression
function lies in the Sobolev space $H^k[0,1]$ (see~\cite{vdVaart07}
for definition), and define the nuisance prior through the
Gaussian process
%
%e7.2 #&#
%
\begin{equation}
\label{eqkIBM}
\eta(t) = \sum_{i=0}^k Z_i \frac{t^i}{i!} + (I_{0+}^k W)(t),
\end{equation}
where $W=\{W_t\dvtx t\in[0,1]\}$ is Brownian motion on $[0,1]$,
$(Z_0,\ldots,Z_k)$ form a~$W$-independent, $N(0,1)$-i.i.d. sample
and $I_{0+}^k$ denotes
$(I_{0+}^1f)(t) = \int_0^tf(s) \,ds$ or $I_{0+}^{i+1}f = I_{0+}^1 I_{0+}^{i}f$
for all $i\geq1$. The prior process $\eta$ is zero-mean Gaussian
of (H\"older-)smoothness $k+1/2$ and the resulting posterior mean for
$\eta$ concentrates asymptotically on the smoothing spline that solves
the penalized ML problem \mbox{\cite{Wahba78,Shen02}}. MCMC simulations based
on Gaussian priors have been carried out by Shively, Kohn and Wood
\cite{Shively99}.

Here, we reiterate the question of how frequentist sufficient conditions
are expressed in a Bayesian analysis based on
Corollary~\ref{corsimplesbvm}. We show that with a nuisance of known
(H\"older-)smoothness greater than $1/2$, the process~(\ref{eqkIBM})
provides a prior such that the marginal posterior for ${\theta}$ satisfies
the Bernstein--von Mises limit.
To facilitate the analysis, we think of the regression function
and the process~(\ref{eqkIBM}) as elements of the Banach space
$(C[0,1],\mbox{$\|\cdot\|_\infty$})$. At a~later stage, we relate
to Banach subspaces with stronger norms to complete the argument.
%
%th7.1 #&#
%
\begin{Thm}
\label{thmplm}
Let $X_1, X_2, \ldots$ be an i.i.d. sample from the partial linear
model~(\ref{eqplrmodel}) with $P_0=P_{{\theta}_0,\eta_0}$ for some
${\theta}_0\in\Tht$, $\eta_0\in H$. Assume\vadjust{\goodbreak} that $H$ is a~subset of
$C[0,1]$ of finite metric entropy with respect to the uniform
norm and that $H$ forms a $P_0$-Donsker class. Regarding the
distribution of $(U,V)$, suppose that $PU=0$, $PU^2=1$ and
$PU^4<\infty$, as well as $P(U-\mathrm{E}[U|V])^2>0$,
$P(U-\mathrm{E}[U|V])^4<\infty$ and
$v\mapsto\mathrm{E}[U|V=v]\in H$. Endow $\Tht$ with a prior that
is thick at ${\theta}_0$ and $C[0,1]$ with a prior $\Pi_H$ such
that $H\subset\operatorname{supp}(\Pi_H)$. Then the
marginal posterior for ${\theta}$ satisfies the Bernstein--von Mises
limit,
%
%e7.3 #&#
%
\begin{equation}
\label{eqplm}
\sup_{B\in\scrB}\bigl|
\Pi\bigl( \sqrt{n}({\theta}-{\theta}_0)\in B\bigmv X_1,\ldots
,X_n \bigr)
- N_{\effDelta_n,\effFI_{{\theta}_0,f_0}^{-1}}(B) \bigr| \convprob
{P_0} 0,
\end{equation}
where $\effscore_{{\theta}_0,\eta_0}(X)=e(U-\mathrm{E}[U|V])$ and
$\effFI_{{\theta}_0,\eta_0}=P(U-\mathrm{E}[U|V])^2$.
\end{Thm}
\begin{pf}
For any ${\theta}$ and $\eta$,
$-P_{{\theta}_0,\eta_0} \log(p_{{\theta},\eta}/p_{{\theta}_0,\eta
_0})=\ft12
P_{{\theta}_0,\eta_0}( ({\theta}-{\theta}_0)U + (\eta-\eta_0)(V))^2$,
so that for fixed ${\theta}$, minimal KL-divergence over $H$
obtains at $\eta^*({\theta}) = \eta_0 - ({\theta}-{\theta}_0)
\mathrm{E}[U|V]$,
$P$-almost-surely. For fixed $\zeta$, the submodel
${\theta}\mapsto Q_{{\theta},\zeta}$ satisfies
%
%e7.4 #&#
%
\begin{eqnarray}
\label{eqparallellik}
&&\log\prod_{i=1}^n\frac{p_{{\theta}_0+n^{-1/2}h_n,\eta^*({\theta}
_0+n^{-1/2}h_n)+\zeta}}
{p_{{\theta}_0,\eta_0+\zeta}}(X_i)\nonumber\\[-3pt]
&&\qquad= \frac{h_n}{\sqrt{n}}\sum_{i=1}^n g_\zeta(X_i)
- \ft12 {h_n}^2 P_{{\theta}_0,\eta_0+\zeta} {g_\zeta}^2\\[-3pt]
&&\qquad\quad{} + \ft12 {h_n}^2 (\PP_n-P)(U-\mathrm{E}[U|V])^2
\nonumber
\end{eqnarray}
for all stochastic $(h_n)$, with $g_\zeta(X)=e(U-{\mathrm E}[U|V])$,
$e=Y-{\theta}_0U-(\eta_0+\zeta)(V)\sim N(0,1)$ under $P_{{\theta
}_0,\eta_0+\zeta}$.
Since $PU^2<\infty$, the last term on the right is
$o_{P_{{\theta}_0,\eta_0+\zeta}}(1)$ if $(h_n)$ is bounded in
probability. We conclude that ${\theta}\mapsto Q_{{\theta},\zeta}$
is stochastically LAN. In addition,~(\ref{eqparallellik}) shows
that $h\mapsto s_n(h)$ is continuous for every $n\geq1$. By assumption,
$\effFI_{{\theta}_0,\eta_0}=P_0{g_0}^2=P(U-{\mathrm E}[U|V])^2$
is strictly positive. We also observe at this stage that $H$ is
totally bounded in $C[0,1]$, so that there exists a constant $D>0$
such that $\|H\|_{\infty}\leq D$.

For any $x\in\RR^3$ and all $\zeta$, the map
${\theta}\mapsto\log{q_{{\theta},\zeta}/q_{{\theta}_0,\zeta
}}(x)$ is continuously
differentiable on all of $\Tht$, with score $g_{{\theta},\zeta}(X)=
e(U-{\mathrm E}[U|V])+({\theta}-{\theta}_0)(U-{\mathrm E}[U|V])^2$. Since
$Q_{{\theta},\zeta}g_{{\theta},\zeta}^2=P(U-{\mathrm E}[U|V])^2
+({\theta}-{\theta}_0)^2P(U-{\mathrm E}[U|V])^4$ does not depend on
$\zeta$
and is bounded over ${\theta}\in[{\theta}_0-\rho,{\theta}_0+\rho]$,
Lemma~\ref{lemUdom} says that $U(\rho_n,h_n)=O(1)$ for all
$\rho_n\downarrow0$ and all bounded, stochastic $(h_n)$. So
for this model, we can apply the rate-free version of the
semiparametric Bernstein--von Mises theorem,
Corollary~\ref{corsimplesbvm}, and its condition (iii)
is satisfied.

Regarding condition (ii) of Corollary~\ref{corsimplesbvm},
we first note that, for $M>0$, $n\geq1$, $\eta\in H$,
\begin{eqnarray*}
\sup_{\|h\|\leq M}-\log\frac{p_{{\theta}_n(h),\eta}}{p_{{\theta
}_0,\eta_0}}
&=& \frac{M^2}{2n}U^2+\frac{M}{\sqrt{n}}\bigl|U\bigl(e-(\eta-\eta
_0)(V)\bigr)
\bigr|\\[-3pt]
&&{}-e(\eta-\eta_0)(V) +\ft12 (\eta-\eta_0)^2(V),\vadjust{\goodbreak}
\end{eqnarray*}
where $e\sim N(0,1)$ under $P_{{\theta}_0,\eta_0}$. With the help of the
boundedness of $H$, the independence of $e$ and
$(U,V)$ and the assumptions on the distribution of $(U,V)$, it is
then verified that condition (ii) of Corollary~\ref{corsimplesbvm}
holds. Turning to condition (i), it is noted that for all
$\eta_1,\eta_2\in H$, $d_H(\eta_1,\eta_2)\leq-P_{{\theta}_0,\eta
_2}\log
(p_{{\theta}_0,\eta_1}/p_{{\theta}_0,\eta_2})=\ft12\|\eta_1-\eta
_2\|_{2,P}^2
\leq\ft12\|\eta_1-\eta_2\|_{\infty}^2$. Hence, for any $\rho>0$,
$N(\rho,\scrP_{{\theta}_0},d_H)\leq N((2\rho)^{1/2},H$, $
\mbox{$\|\cdot\|_{\infty}$})<\infty$. Similarly, one shows that for
all $\eta$ both $-P_0\log(p_{{\theta}_0,\eta}/p_{{\theta}_0,\eta
_0})$ and
$P_0(\log(p_{{\theta}_0,\eta}/p_{{\theta}_0,\eta_0}))^2$ are
bounded by
$(\ft12+D^2)\|\eta-\eta_0\|_{\infty}^2$. Hence, for any
$\rho>0$, $K(\rho)$ contains a \mbox{$\|\cdot\|_{\infty}$}-ball. Since
$\eta_0\in\operatorname{supp}({\Pi_H})$, we see that condition (i) of
Corollary~\ref{corsimplesbvm} holds. Noting that
$({p_{{\theta}_n(h),\eta}}/{p_{{\theta}_0,\eta}}(X))^{1/2}=
\exp((h/2\sqrt{n})eU-(h^2/ 4n)U^2)$, one derives
the $\eta$-independent upper bound,
\[
H^2\bigl( P_{{\theta}_n(h_n),\eta} , P_{{\theta}_0,\eta} \bigr)
\leq\frac{M^2}{2n}PU^2 + \frac{M^3}{6n^2} PU^4 = O(n^{-1})
\]
for all bounded, stochastic $(h_n)$, so that condition (iv)
of Corollary~\ref{corsimplesbvm} holds.

Concerning condition (v), let $(M_n)$,
$M_n\rightarrow\infty$ be given and define $\Tht_n$ as in
Section~\ref{secmarg}. Rewrite $\sup_{\eta\in H}
\sup_{{\theta}\in\Tht_n^c}\PP_n\log(p_{{\theta},\eta
}/p_{{\theta}_0,\eta})
= \sup_{{\theta}\in\Tht_n^c}( ({\theta}-{\theta}_0)\times(\sup_\zeta
\PP_nZW)
- \ft12 ({\theta}-{\theta}_0)^2 \PP_nW^2)$,
where $Z=e_0-\zeta(V)$, $W=U-{\mathrm E}[U|V]$. The maximum-likelihood
estimate $\hat{{\theta}}_n$ for ${\theta}$ is therefore of the form
$\hat{{\theta}}_n={\theta}_0 + R_n$, where $R_n=\sup_\zeta\PP
_nZW / \PP_nW^2$.
Note that $P_0ZW=0$ and that $H$ is assumed to be $P_0$-Donsker,
so that $\sup_\zeta\GG_nZW$ is asymptotically tight. Since, in addition,
$\PP_nW^2\rightarrow P_0W^2$ almost surely and the limit is strictly
positive by assumption, $P_0^n( \sqrt{n} |R_n|> \ft14M_n) = o(1)$.
Hence,
\begin{eqnarray*}
&&
P_0^n\biggl( \sup_{\eta\in H}\sup_{{\theta}\in\Tht_n^c}
\PP_n\log\frac{p_{{\theta},\eta}}{p_{{\theta}_0,\eta}} > -\frac
{CM_n^2}{n}\biggr)\\[-2pt]
&&\qquad\leq P_0^n\biggl( \sup_{{\theta}\in\Tht_n^c}
\biggl( \ft14|{\theta}-{\theta}_0|\frac{M_n}{n^{1/2}} - \ft12
({\theta}-{\theta}
_0)^2\biggr)
\PP_nW^2 > -\frac{CM_n^2}{n} \biggr) + o(1)\\[-2pt]
&&\qquad\leq P_0^n( \PP_nW^2<4C )+o(1).
\end{eqnarray*}
Since $P_0W^2>0$, there exists a $C>0$ small enough such that the first
term on the right-hand side is of order $o(1)$ as well, which shows that
condition~(\ref{eqlehmann}) is satisfied. Lemma~\ref{lemlehmann}
asserts that condition (v) of Corollary~\ref{corsimplesbvm}
is met as well. Assertion~\ref{eqplm} now holds.\vspace*{-2pt}
\end{pf}

In the following corollary we choose a prior by picking a suitable
$k$ in~(\ref{eqkIBM}) and conditioning on $\|\eta\|_{\alpha}<M$.
The resulting prior is shown to be well defined below and is denoted
$\Pi^k_{\alpha,M}$.\vspace*{-2pt}
%
%co7.2 #&#
%
\begin{corollary}
\label{corsmoothplr}
Let $\alpha>1/2$ and $M>0$ be given; choose
$H=\{\eta\in C^\alpha[0,1]\dvtx\|\eta\|_\alpha<M\}$ and assume that
$\eta_0\in C^\alpha[0,1]$. Suppose the distribution of the covariates
$(U,V)$ is as in Theorem~\ref{thmplm}. Then, for any integer
$k>\alpha-1/2$, the conditioned prior $\Pi^k_{\alpha,M}$ is
well defined and gives rise to a marginal posterior for ${\theta}$
satisfying~(\ref{eqplm}).
\end{corollary}
\begin{pf}
Choose $k$ as indicated; the Gaussian distribution of $\eta$ over
$C[0,1]$ is based on the RKHS $H^{k+1}[0,1]$ and denoted $\Pi^k$.
Since $\eta$ in~(\ref{eqkIBM}) has smoothness $k+1/2>\alpha$,
$\Pi^k(\eta\in C^\alpha[0,1])=1$. Hence, one may also
view $\eta$ as a Gaussian element in the H\"older class
$C^\alpha[0,1]$, which forms a separable Banach space even with
strengthened norm $\mbox{$\|\cdot\|$}=\|\eta\|_{\infty}+
\mbox{$\|\cdot\|_{\alpha}$}$,
without changing the RKHS. The trivial
embedding of $C^\alpha[0,1]$ into $C[0,1]$ is one-to-one and
continuous, enabling
identification of the prior induced by $\eta$ on $C^\alpha[0,1]$
with the prior $\Pi^k$ on $C[0,1]$. Given $\eta_0\in C^\alpha[0,1]$
and a sufficiently smooth kernel $\phi_\sigma$ with bandwidth
$\sigma>0$, consider $\phi_\sigma\star\eta_0\in H^{k+1}[0,1]$.
Since $\|\eta_0-\phi_\sigma\star\eta_0\|_{\infty}$ is of order
$\sigma^\alpha$, and a similar bound exists for the $\alpha$-norm
of the difference~\cite{vdVaart07}, $\eta_0$ lies in the closure
of the RKHS both with respect to \mbox{$\|\cdot\|_{\infty}$} and to
\mbox{$\|\cdot\|$}. Particularly, $\eta_0$ lies in the support of $\Pi^k$,
in $C^\alpha[0,1]$ with norm \mbox{$\|\cdot\|$}. Hence,
\mbox{$\|\cdot\|$}-balls centered on $\eta_0$ receive nonzero prior mass,
that is, $\Pi^k(\|\eta-\eta_0\|<\rho)>0$ for all $\rho>0$. Therefore,
$\Pi^k(\|\eta-\eta_0\|_\infty<\rho,
\|\eta\|_\alpha<\|\eta_0\|_\alpha+\rho)>0$, which guarantees
that $\Pi^k(\|\eta-\eta_0\|_\infty<\rho,\|\eta\|_\alpha<M)>0$,
for small enough $\rho>0$. This implies that $\Pi^k(\|\eta\|_\alpha<M)>0$,
and
\[
\Pi^k_{\alpha,M}(B) = \Pi^k( B\bigmv \|\eta\|_\alpha<M )
\]
is well defined for all Borel-measurable $B\subset C[0,1]$. Moreover,
it follows that $\Pi^k_{\alpha,M}(\|\eta-\eta_0\|_\infty<\rho)>0$ for
all $\rho>0$. We conclude that $k$ times integrated Brownian motion
started at random, conditioned to be bounded by $M$ in $\alpha$-norm,
gives rise to a prior that satisfies
$\operatorname{supp}(\Pi^k_{\alpha,M})=H$. As is well-known~\cite{vdVaart96},
the entropy numbers of $H$ with respect to the
uniform norm satisfy, for every $\rho>0$, $N(\rho,H,\mbox{$\|\cdot\|
_{\infty}$}) \leq K\rho^{-1/\alpha}$, for some constant $K>0$ that
depends only on $\alpha$ and~$M$. The associated bound on the
bracketing entropy gives rise to finite bracketing integrals, so that
$H$ universally Donsker. Then, if the distribution of the covariates
$(U,V)$ is as assumed in Theorem~\ref{thmplm}, the Bernstein--von Mises
limit~(\ref{eqplm}) holds.
\end{pf}

\section*{Acknowledgments}

The authors would like to thank D. Freedman,\break A.~Gamst, C. Klaassen, B.
Knapik and A. van der Vaart for valuable discussions and suggestions.
B. J. K. Kleijn thanks U.C. Berkeley's Statistics Dept. and Cambridge's
Isaac Newton Institute for their hospitality.

%suskaldyti doi

% imsref loaded by lrinkeviciute, 2012-01-19 10:19:37
% imsref loaded by lrinkeviciute, 2012-01-19 10:43:46
%

\printaddresses

\end{document}